\documentclass[12pt]{amsart}
\usepackage{amsmath,amssymb,slashed}
\usepackage{latexsym}
\usepackage{delarray}
\usepackage{mathrsfs}
\usepackage{color}
\usepackage{dsfont,float}

\usepackage{hyperref}
\usepackage{graphicx}

\numberwithin{equation}{section}
\theoremstyle{plain}
\newtheorem{thm}{Theorem}[section]
\newtheorem{theorem}{Theorem}
\newtheorem{lem}[thm]{Lemma}
\newtheorem{prop}[thm]{Proposition}
\newtheorem{cor}[thm]{Corollary}
\newtheorem{exam}[thm]{Example}
\theoremstyle{definition}
\newtheorem{defn}{Definition}[section]
\newtheorem{remark}{Remark}[section]

\newcommand{\vertiii}[1]{{\left\vert\kern-0.25ex\left\vert\kern-0.25ex\left\vert #1
    \right\vert\kern-0.25ex\right\vert\kern-0.25ex\right\vert}}

\newcommand\R{{\mathbb R}}
\newcommand\C{{\mathbb C}}

\newcommand\Rn{{{\mathbb R}^n}}
\newcommand\ep{\varepsilon}
\newcommand\va{\varphi}
\newcommand\pa{\partial}
\def\im{{\rm i}}

\title[Local energy decay estimate]{On wave equation outside trapping obstacles and local energy decay in odd space dimensions}
\author[Vladimir Georgiev]{Vladimir Georgiev}
\address{
   Vladimir Georgiev:
   \endgraf
  Dipartimento di Mathematica, Universit\`a di Pisa
  \endgraf
  Largo Bruno Pontecorvo 5, 56127
  \endgraf Pisa, Italia
  \endgraf
 and
  Faculty of Science and Engineering, Waseda University
   \endgraf
3-4-1, Okubo, Shinjuku-ku
\endgraf
Tokyo 169-8555
\endgraf
Japan
 \endgraf
and IMI--BAS, Acad.
Georgi Bonchev Str., Block 8, 1113
\endgraf
Sofia, Bulgaria
\endgraf
  {\it E-mail address} {\rm georgiev@dm.unipi.it}
  }

\author[Tokio Matsuyama]{Tokio Matsuyama}
\address{
Tokio Matsuyama:
 \endgraf
Department of Mathematics
\endgraf
Chuo University
\endgraf
1-13-27, Kasuga, Bunkyo-ku
\endgraf
Tokyo 112-8551
\endgraf
Japan
\endgraf
{\it E-mail address} {\rm tokio@math.chuo-u.ac.jp}
}

\thanks
{
\endgraf  
The first author was supported by Gruppo Nazionale per l'Analisi Matematica 2020, by the project PRIN 2020XB3EFL with the Italian Ministry of Universities and Research, by Institute of Mathematics and Informatics, Bulgarian Academy of Sciences, by Top Global University Project, Waseda University and the Project PRA 2022 85 of University of Pisa.
The second author was supported by Grant-in-Aid for Scientific Research (C) (No. 18K03377),
Japan Society for the Promotion of Science.}


\begin{document}

\begin{abstract}
The purpose of the present paper is to establish appropriate
cut-off resolvent estimates for the Dirichlet Laplacian on exterior domains.
The  geometrical assumptions on domains are rather general,
for example, non-trapping condition is not imposed.
The first key assumption guarantees the result on propagation of singularities, 
the second is the smallness of the Lebesgue measure of the portions of the trapped sets in the fibers of cosphere bundle,
and the third concerns the upper bound of the sojourn time. 
As a by-product of these estimates, the local energy decay estimate for
solutions to the initial-boundary value problem for wave equation in the case of odd space dimensions is obtained.
\end{abstract}

\maketitle


\section{Introduction} \label{sec:1}
Let $\Omega$ be an open connected set in $\R^n$ which is the exterior of a compact obstacle
\[
 \mathcal{O}= \R^n \setminus \Omega
\]nifty
with a $C^\infty$-boundary $\pa \Omega$.
We denote by $ \left. \Delta\right\vert_D$  the standard Laplace operator  with the homogeneous Dirichlet boundary condition in $\Omega .$
Let us consider the initial-boundary value problem for wave equation
of the form
\begin{equation}\label{eq.1}
\left\{
\begin{aligned}
    & \partial_{t}^2 u - \left. \Delta\right\vert_D u = 0,  &\quad (t,x)\in \R\times \Omega,\\
   & u(0,x)=u_0(x),\quad  \partial_t u(0,x) = u_1(x), &\quad x\in \Omega.
\end{aligned}
\right.
\end{equation}
Our first goal is to obtain the estimates for the following weighted resolvent operator
\begin{equation} \label{EQ:Resolvent}
\langle x \rangle^{-s} \left(-\left. \Delta\right\vert_D - z^2\right)^{-1}  \langle x \rangle^{-s}, \quad s>1/2
\end{equation}
for $\mathrm{Re}(z)\gg1$ and $|\mathrm{Im}(z)|\ll1$, where we put
\[
\langle x\rangle=\sqrt{1+|x|^2}.
\]
It is well known that when
$\Omega$ is a non-trapping domain, $L^2$-$L^2$-estimate for the operator \eqref{EQ:Resolvent} has a bound
like $1/\mathrm{Re}(z)$ for $\mathrm{Re}(z)>0$ (see, e.g.,
Mochizuki \cite{Mochizuki-2,Mochizuki-1}).
When $\Omega$ is generic (possibly  trapping) domain, Burq in \cite{Burq-Acta} proved that the cut-off resolvent
\begin{equation} \label{w-resolvent}
 \chi \left(-\left. \Delta\right\vert_D - z^2\right)^{-1}  \chi, \quad \chi \in C_0^\infty
 \end{equation}
has an exponential bound like $e^{c\mathrm{Re}(z)}$ for $\mathrm{Re}(z)>0$ (see \cite{Burq-Acta}).
As a by-product of these estimates, it is proved in Theorem 1 from \cite{Burq-Acta} that the local energy decays in a logarithmic rate. A similar result for general aymptotically flat manifolds is obtained  in \cite{Moschidis}.\\

The aim in this paper is to prove that geometrical requirements imply the uniform resolvent estimate
for the weighted resolvent operator \eqref{EQ:Resolvent}, and then one can get the uniform bound for the cut-off resolvent operator \eqref{w-resolvent}, when $\Omega$ is a trapping domain.
The proof is based on an appropriate micro-localization of
the source term $f$ in the equation
\begin{equation*} 
   \left( \left. -\Delta \right\vert_D - (\lambda^2 \pm \im \varepsilon) \right)w =f.
\end{equation*}
If singularities associated with $f$ are outside the trapped set (their distance to the trapped set is bounded from below by $\delta>0$), then we use an assumption giving the nice bound of the sojourn time, and via propagation of singularity analysis we get smoothing estimates for $w.$
If the localized portion of $f$ has singularity in $\delta$-neighbourhood of the trapped set, then we use the exponential
growth estimate due to \cite{Burq-Acta} and the smallness of the Lebesgue measure of the portions of the trapped sets in the fibers of cosphere bundle.\\

Once the resolvent bound are obtained, one can deduce the local energy decay estimates for
solutions to the initial-boundary value problem for wave equation.
The local energy for wave equations
is defined by letting
\[
E_R(t) = \int_{\Omega_R} \left[ |\nabla u(t,x)|^2
 + |\partial_t u(t,x)|^2 \right] dx,
\]
where we set
\[
\Omega_R:=\Omega \cap \{ |x| \leq R\}.
\]
Here and below, $R>0$ is chosen so that
\[
 \mathcal{O} = \R^n \setminus \Omega \subseteq \{ |x| \leq R \}.
 \]
The result due to Ralston \cite{R} concerns the case that
\[
\text{$\mathcal{O}$ is a compact and trapping obstacle in $\R^3$,}
\]
and his result asserts that, given any $ \mu
\in(0,1)$ and $T > 0,$ one can find initial data
 $u(0), \partial_t u(0) \in C_0^\infty{(\Omega)}$ with
\[ \int_{\Omega }\left[ |\nabla u(0,x)|^2 + |\partial_t u(0,x)|^2 \right] dx = 1
\]
such that the solution to the initial-boundary value problem \eqref{eq.1}
satisfies the inequality
\[
 E_R(T) \geq 1 -\mu.
 \]
Ralston's result  is rewritten also in the following way:
\begin{equation*}
 \sup_{t>0, \ E_R(0) = 1} E_R(t)  \geq 1-\mu,
\end{equation*}
provided that $(u(0,\cdot), \partial_t u(0,\cdot))$ are compactly supported in $\Omega_R.$
The  results of Burq \cite{Burq-Acta} and Ralston \cite{R} show that we really need better regularity for the initial data.
For this purpose, we introduce the local energy norm of initial data with $\alpha$ loss of regularity
\[
 E^{(\alpha)}_R (0): = \|u(0)\|^2_{H^{\alpha+1}(\Omega_R)} + \|\partial_tu(0)\|^2_{H^\alpha(\Omega_R)}
 \]
and
we can rewrite the result in \cite{Burq-Acta} in the form
\begin{equation*}
 \sup_{t\ge 0, \, E^{(\alpha)}_R(0) \leq 1 } \log (2+t)^{\alpha} \cdot E_R(t)   < \infty,
\end{equation*}
for $\alpha >0$, and again $(u(0,\cdot), \partial_t u(0,\cdot))$ are compactly supported in $\Omega_R.$
Our geometrical assumptions (see Assumptions A, B and C in Section \ref{sec:2}) guarantee that if $n$ is an odd integer
with $n\geq 3$, then there exists $\kappa>0$ such that 
\begin{equation*}
 \sup_{t\ge0,\, E^{(\alpha)}_R(u)(0) \leq 1 } e^{\kappa t} E_R(t)  < \infty
\end{equation*}
for $\alpha >0$, and again $(u(0,\cdot), \partial_tu(0,\cdot))$ are compactly supported in $\Omega_R.$

On the other hand, the scattering theory developed by Lax and Phillips gives a
construction of the scattering operator by using weaker form of local energy decay\begin{equation}\label{eq.LED2}
   \liminf_{t \to \infty} E_R(t) =0
 \end{equation}
(see \cite{LP}, and also  \cite{Petkov}).
Note that \eqref{eq.LED2} follows directly from the RAGE (or simply ergodic type) theorem
\begin{equation} \label{eq.LED2a}
    \lim_{T \to \infty} \frac{1}{T} \int_0^T E_R(t) dt = 0
\end{equation}
and the property that zero is not eigenvalue of
$\left. \Delta\right\vert_D$, i.e.,
 \[
   u \in \mathrm{Dom}(\left. \Delta\right\vert_D), \, \left. \Delta\right\vert_D u=0 \quad  \Longrightarrow \quad u=0.
\]
An important consequence of weak energy decay \eqref{eq.LED2a} is the existence of the wave operators
\[
  W_\mp: = \mathrm{s-}\lim_{t \to \pm \infty}    e^{\im t \sqrt{-\left. \Delta\right\vert_D}} J_0 e^{-\im t \sqrt{-\Delta}},
\]
where $J_0$ is  the orthogonal projection
\[
 J_0 : L^2(\R^n) \longrightarrow L^2(\Omega).
\]
This observation implies that scattering theory and existence of wave operators are established without appealing to the additional geometric assumption of type that
\begin{equation}\label{eq.LED5}
   \mbox{$\mathcal{O} = \R^n \setminus \Omega$ is a non-trapping obstacle}.
\end{equation}
The condition \eqref{eq.LED5} is crucial for the strong local energy decay in view of  the results of Morawetz, Ralston and Strauss \cite{MRS} and Ralston \cite{R}.\\

Our main decay estimates below
are obtained also without appealing to
assumption \eqref{eq.LED5} and these are
probably the main novelty in
our work.\\

This paper is organized as follows:
In Section \ref{sec:2}, after recalling notions arising in symplectic geometry, and giving
a definition of trapped sets, we state main results. In Section \ref{sec:3}
we prepare some results on propagation of singularities of solutions to wave equation.
Section \ref{sec:4} is devoted to deriving resolvent estimates in high frequency.
In Section \ref{sec:5} we prove the resolvent estimates.
Section \ref{sec:6} is devoted to proving the proof of local energy decay estimate.
The resolvent estimate in low frequency is proved in Section \ref{sec:7}. 
In Appendix some examples of trapped domains fulfilling our assumptions are provided, 
and we prove the commutator estimates and recall some basic resolvent estimates in the free case.


\section{Main results} \label{sec:2}

Let us introduce the local behaviour of the Hamiltonian flow
and the notion of the trapped set.
We start by recalling some basic notions from symplectic geometry.
We denote by $\overline{\Omega}$ the closure of domain $\Omega,$ i.e.,
\[
\overline{\Omega} = \Omega \cup \partial \Omega.
\]
Let $T^*(\overline{\Omega})$ be the cotangent bundle of $\overline{\Omega}$
and denote by $S^*(\overline{\Omega})$ the corresponding cosphere bundle (see Melrose and Sj\"ostrand \cite{MS78}, and also Section 24.3
in H\"ormander \cite{Hormander}).
Namely, if $$\rho = (x, \xi) \in  T^*(\overline{\Omega}),$$
then
\[
\rho = (x, \xi) \in \overline{\Omega} \times  (\R^n\setminus \{0\}),
\]
and moreover,
\[
\rho=(x,\xi) \in S^*(\overline{\Omega})
\]
means that
\[
\rho=(x,\xi) \in \overline{\Omega}\times \mathbb{S}^{n-1}.
\]
Then we can define a canonical metric $d$ on  $T^*(\overline{\Omega})$ induced by the standard Euclidean metric on $\Omega.$
The corresponding Hamiltonian flow
\[
 \mathcal{H}_p^t : \rho \in S^*(\overline{\Omega}) \longrightarrow \mathcal{H}_p^t(\rho) \in
S^*(\overline{\Omega})
\]
associated with \eqref{eq.1} is defined by the aid of the symbol $p(\tau,\xi) = \tau^2 - |\xi|^2$ of
the operator $\pa^2_t-\left. \Delta\right\vert_D$. By using the standard projection operator
\[
 \Pi: T^*(\overline{\Omega}) \longrightarrow \overline{\Omega},
 \]
it is easy to see that the Hamiltonian orbit
$ \{ \mathcal{H}^t_p(\rho); t \in \R \}  $ has the projection $ \Pi \,(\mathcal{H}_p^s(\rho))$ containing the union of segments in $\Omega.$
 It should be mentioned that Melrose and Sj\"ostrand established the existence of broken Hamiltonian flow associated to  $p(\tau,\xi)$ (see Section 3 in \cite{MS78}).
The Hamiltonian flow is related to the propagation of singularities of the solution due to \cite{MS78,MS82}.\\

Our first assumption guarantees the result on the propagation of singularities due to Theorem (0.11) from Melrose and
Sj\"ostrand (see \cite{MS78}, and also \cite{MS82}).\\

\noindent
{\bf Assumption A.}
The Hamiltonian flow associated with  $p(\tau,\xi)= \tau^2-|\xi|^2,$ starting from non-trapping points is tangent to
$\pa T^*(\Omega)\setminus \{0\}$ at $p=0,$ to finite order. \\

We shall make further two key assumptions.
The trapped set is formed by Hamiltonian orbits having compact projection on $x$-space.
More precisely, we introduce the following:
\begin{defn}
$\rho \in S^*(\overline{\Omega})$ is said to belong to the {\em future trapped set} $\Sigma_+$ if
there exist a compact set $K_+$ in $\overline{\Omega}$ and an instant $t^+_0\ge0$ such that
\[
\Pi (\mathcal{H}^t_p(\rho))\subseteq K_+
\]
for all $t\ge t^+_0$.
Similarly, $\rho \in S^*(\overline{\Omega})$ is said to belong to the {\em past trapped set} $\Sigma_-$ if
there exist a compact set $K_-$ in $\overline{\Omega}$ and an instant
$t^-_0\le0$ such that
\[
\Pi (\mathcal{H}^t_p(\rho))\subseteq K_-
\]
for all $t\le t^-_0$. The trapped set $\Sigma$ is defined by
\[
\Sigma = \Sigma_+ \cup \Sigma_-.
\]
\end{defn}

\begin{remark}
It is clear that $\rho \in \Sigma_\pm$ if and only if the  projection
$ \Pi\left( \left\{ \mathcal{H}_p^t(\rho)  ; t \gtrless 0 \right\} \right)$
is compact
in $\overline{\Omega}.$
Moreover, $\Sigma_+$ is invariant under the action of $\mathcal{H}^t_p$ with respect to
$t \geq 0,$ i.e., $\mathcal{H}^t_p(\Sigma_+) \subseteq \Sigma_+$ for any $t \geq 0.$
\end{remark}

We observe that for any $\delta >0$, $R\geq 1$ and $\rho \in S^*(\overline{\Omega})$ outside the trapped set $\Sigma$ with
$ d(\rho, \Sigma) \geq \delta ,$ there exists $t_{\delta, R, \rho}>0$ such that
\begin{align*}
\text{(i)} & \quad \text{$\mathrm{dist} \left( \Pi (\mathcal{H}^{t_{\delta, R, \rho}}_p(\rho)) , \Pi(\Sigma)\right) \geq  R$ \quad  for $t \geq t_{\delta,R,\rho};$}  \\
\text{(ii)} & \quad \text{$\mathrm{dist}\left( \Pi (\mathcal{H}^{t_{\delta, R, \rho}}_p(\rho)) , \Pi(\Sigma)\right) <  R$ \quad for $0\le t<t_{\delta,R,\rho}$,}
\end{align*}
i.e., $t_{\delta,R,\rho}$ is the first moment when the property (i) is fulfilled.
Here $\mathrm{dist}(\cdot,\cdot)$ is the standard Euclidean distance on $\Omega$.

\begin{defn}
$t_{\delta,R,\rho}$ is said to be a {\em sojourn time} if it satisfies {\rm (i)} and {\rm (ii)}.
\end{defn}

The second assumption guarantees that a small neighborhood of $\Sigma$ has a small symplectic measure.
We recall that $T^*(\overline{\Omega})$ is the bundle with projection $$\Pi: T^*(\overline{\Omega}) \longrightarrow  \overline{\Omega} $$ so that for any $\delta$-neighbourhood $\Sigma_\delta$ of the trapped set $\Sigma$, we  can evaluate the Lebesgue measure of $\Sigma_\delta \cap \Pi^{-1}(x)$ in the fiber $\Pi^{-1}(x)$ for any $x \in \Pi(\Sigma_\delta)$.  

\vspace{5mm} 

The following assumption gives more precise information about $\Sigma_\delta.$\\

\noindent 
{\bf Assumption B.} There exists $\delta_0>0$ such that for any $\delta \in (0,\delta_0)$, the corresponding 
$\delta$-neighbourhood $\Sigma_\delta$ of the trapped set $\Sigma$ fulfills the following two properties:
\begin{enumerate}
\item[(i)] 
$\Sigma_\delta$ is represented as  
\begin{equation} \label{eq.cove}
 \Sigma_\delta = \bigcup_{j=1}^{N(\delta)} (K_{j,\delta} \times E_{j,\delta}),
 \end{equation}
where $K_{j,\delta}$ are compact subsets of $\Omega$, 
and $E_{j,\delta}$ are  balls with radii proportional to $\delta$ in
$\mathbb{S}^{n-1}.$
\item[(ii)] 
 There exist constants $C>0$ and $\nu \in [0,n-1)$ 
such that  
\begin{equation} \label{eq.impcm1}
  N(\delta) \leq C \delta^{-\nu}.
 \end{equation}
\end{enumerate}
\vspace{5mm} 
\begin{remark}
We can interpret Assumption B by using the notion of $\delta$-external covering of the fibre projection of the set 
$\Sigma_\delta.$
    More precisely,  
    for the finite set $\{1,2,\cdots, N(\delta)\}$ we can assert that for any $x \in \Pi(\Sigma_\delta)$ we have the covering of the fiber over $x$
    \begin{equation*} 
\Pi^{-1} (x) = \bigcup_{j=1}^{N(\delta)}\left( \{x\} \times E_{j,\delta} \right) ,
 \end{equation*}
   and therefore, we have $\delta$-external covering of cardinality $N(\delta)$ satisfying \eqref{eq.impcm1}.
    \end{remark}

In many cases, it is checked that $\Sigma_\delta$ is a trivial bundle, and Assumption B is reduced to the following: \\ 

\noindent
{\bf Assumption B${}^\prime$.} There exists $\delta_0>0$ such that for any $\delta \in (0,\delta_0)$ the corresponding $\delta$ neighbourhood $\Sigma_\delta$ of the trapped set $\Sigma$ fulfills the following two properties:
\begin{enumerate}
\item[(i)] 
$\Sigma_\delta$ is trivial, i.e.
$$ \Sigma_\delta = K_\delta \times E_{\delta},$$
where $K_{\delta}$ is a compact subset of $\Omega$, 
and $E_{\delta}$ is a subset of ball with radius proportional to $\delta$ in
$\mathbb{S}^{n-1}.$

\item[(ii)] There exist constants $C>0$ and $\nu \in [0,n-1)$
such that  the Lebesgue measure  $\mu(E_\delta)$ of $E_\delta$ satisfies
$$ \mu(E_\delta) \leq C \delta^{n-1-\nu}.$$
\end{enumerate}

\vspace{2mm}

\begin{remark}
Clearly, Assumption B${}^\prime$ implies Assumption B, and in the case of two convex obstacles in 3-dimensional spaces we see that 
$\mu(E_\delta) \sim \delta^2$. Hence, in this case we have Assumption B${}^\prime$ with $\nu =0.$
\end{remark}

Let  $R>0$ be such that the ball $\{x\in\Rn:|x|<R\}$ covers the obstacle $\mathcal{O}$. 
Our key geometrical assumption concerning the behaviour of the Hamiltonian flow near $\Sigma$ is the following:

\vspace{2mm}

\noindent
{\bf Assumption C.}
There exist  $\delta_1>0$, 
$N_0 \geq 0$ and $C(R)>0$ such that for any $\delta \in (0,\delta_1)$ and $\rho \in S^*(\overline{\Omega})$ with
$d(\rho, \Sigma) \geq \delta$, the sojourn time $t_{\delta,R,\rho}$ satisfies
\begin{equation*} 
 t_{\delta,R,\rho} \leq C(R) \left[\log \left( \frac{1}{\delta}\right)\right]^{N_0}.
\end{equation*}
 
In Appendix \ref{Append:Appendix A} several examples of domains satisfying or not satisfying the above assumptions are provided. \\

We denote by $\mathcal{B}(X,Y)$ the space of all bounded linear operators from a Banach space $X$ into another one $Y$, and
write $\mathcal{B}(X)=\mathcal{B}(X,X)$. $\|T\|_{\mathcal{B}(X,Y)}$ stands for the operator norm of $T\in \mathcal{B}(X,Y)$. \\

Our first main result reads as follows.
\begin{theorem}\label{thm:main}
Let $n$ be an odd integer with $n\ge3$, and $\Omega$ be a domain in $\R^n$ which is
the exterior of a compact obstacle $\mathcal{O}$ with a smooth boundary
$\partial \Omega$. Suppose that $\Omega$ satisfies Assumptions {\rm{A}}, {\rm{B}} and {\rm{C}}.
Then there exists a sufficiently small constant $\varepsilon_0 >0$ such that for any cut-off function $\chi$ which is identically $1$ in a neighbourhood of the obstacle
$\mathcal{O}$, the operator
\[
\chi  \left(-\left. \Delta\right\vert_D - z^2\right)^{-1}  \chi
\]
has no poles in the region
$$\Gamma(\ep_0):=\{z\in \C: z^2=\lambda^2 \pm \im \varepsilon, \, \lambda > 0, \, 0 <\varepsilon 
\leq \varepsilon_0\}.$$
Furthermore, there exists a constant $C>0$ such that
\begin{equation} \label{eq.ME1}
\sum_{k=0}^1 |z|^{1-k} \big\|\chi\nabla^k
\left(-\left. \Delta\right\vert_D - z^2\right)^{-1}  \chi
\big\|_{\mathcal{B}(L^2(\Omega))}\leq C \log(2+\lambda)
\end{equation}
for any $z\in \Gamma(\ep_0)$.
\end{theorem}

The second is as follows:
\begin{theorem}\label{thm:mainB}
Under the assumptions of Theorem {\rm \ref{thm:main}}, there exists a constant $C>0$ such that
\begin{equation} \label{eq.ME1n1}
\sum_{k=0}^1 |z|^{1-k} \big\|\chi\nabla^k
\left(-\left. \Delta\right\vert_D - z^2\right)^{-1}  \chi f 
\big\|_{L^2(\Omega)}\leq C \|f\|_{H^2(\Omega)}
\end{equation}
for any $z\in \Gamma(\ep_0)$ and $f\in H^2(\Omega).$
\end{theorem}

As a consequence of Theorem \ref{thm:mainB}, we have the following:
\begin{theorem} \label{thmLED}
Let $n$ and $\Omega$ be as in
Theorem {\rm \ref{thm:main}}. Consider
\[
\Omega_R = \Omega \cap \{|x| \leq R\},
\]
where $R>0$ is chosen so large that
\[
 \mathcal{O} \subseteq \{ |x| \leq R \}.
 \]
If a pair of functions $(u_0,u_1)$ in \eqref{eq.1} belongs to $[H^{\alpha+1}(\Omega)\cap H^1_0(\Omega)]\times
H^{\alpha}(\Omega)$ for some $\alpha>0$,
and if it is compactly supported in
 $\Omega_R$, then there exist positive constants $C$ and $\kappa$ such that the solution $u$ to the
initial-boundary value problem \eqref{eq.1} satisfies  the estimate
\begin{equation}\label{eq.led}
    E_R(t) \leq C e^{-\kappa t}\left(
    \|u_0\|^2_{H^{\alpha+1}(\Omega)} + \|u_1\|^{2}_{H^{\alpha}(\Omega)}  \right)
 \end{equation}
 for any $t > 0$.
 \end{theorem}

Theorem \ref{thmLED} covers the result of Ikawa (see \cite{Ikawa}). Indeed, if the obstacle $\mathcal{O}$ consists of several disjoint strictly convex bodies in $\R^3$, then it is proved in \cite{Ikawa} that the local energy decays exponentially provided that the initial data $(u_0,u_1)$ are compactly supported and belong to $H^3(\Omega)\times H^2(\Omega)$.
If $\Omega$ is an exterior domain in even dimensional Euclidean spaces with $n\geq 4$,
the local energy decays polynomially. The proof can be found in a forthcoming work.\\

In this paper we often write $f\lesssim g$ or $g\gtrsim f$, if there exists a positive constant $C$ such that $f\leq C g$. If $f\lesssim g$ and $g\lesssim f$, we write $f\sim g$. Furthermore, we denote by $[A,B]:=AB-BA$ the commutator of operators $A$ and $B$.


\section{Propagation of singularities}
\label{sec:3}
In this section we introduce the smoothing estimates for the exterior problem \eqref{eq.1}.
We start by recalling the definition of the wave front set
$\mathrm{WF}_b(u)$ of $u\in \mathscr{D}^\prime(\Omega)$ introduced by H\"ormander (see \cite{Hormander}).

Let $x_0\in\Omega$.
 We say that
 $(x_0,\xi_0)\notin \mathrm{WF}(u)$
if there exist smooth functions $a(x)$ and
real-valued $\va(x,\xi)$ such that
\[
\text{$\mathrm{supp}\, a(\cdot)$ is compact,}
\]
\[
a(x_0)\ne0, \quad \nabla_x \va(x_0,\xi_0)=\xi_0,
\]
\[
\left|\int_{\Omega} e^{\mathrm{i}\va(x,\xi)}a(x)u(x)\, dx\right|\lesssim |\xi|^{-M} \quad
\quad \text{for $\frac{\xi}{|\xi|}$ close to $\frac{\xi_0}{|\xi_0|}$ and arbitrary $M>0$}
\]
(see p.~15 in \cite[Duistermaat]{Duistermaat}).
Now, referring to the definition of $\mathrm{WF}(u)$, we shall state the definition of
$\mathrm{WF}_b(u)$.

\begin{defn}[\cite{Hormander}]
We say that $\rho\notin \mathrm{WF}_b(u)$
if there exists an open neighbourhood $U$ of
$\rho\in T^*(\overline{\Omega})$ with compact closure $\overline{U}$
such that $Pu\in C^\infty(\overline{\Omega})$ for any pseudo-differential operator $P$
with symbol supported in $U$.
\end{defn}

For any compact set $K \subseteq S^*(\overline{\Omega})$, the inclusion
\[
\mathrm{WF}_b(u) \cap S^*(\overline{\Omega})\subseteq K
\]
means that
\[
Pu \in C^\infty(\overline{\Omega})
\]
holds for any pseudo-differential operator $P=P(x,D_x)$ of order zero satisfying
\[ \mathrm{supp} \, p \cap K = \emptyset,
\]
where $p=p(x,\xi)$ is the symbol of $P(x,D_x).$
We remark that the class of zero order pseudo-differential operators with symbols in conic neighborhoods in $T^*(\Omega)$ was introduced and studied by Melrose and  Sj\"ostrand (see \cite{MS78,MS82}, and also Mazzeo and Melrose \cite{MM87}).

To state the result on the propagation of singularities,
we use the notation:
\[
B_R=\{x\in \Rn:|x|<R\}
\]
the  ball with radius $R$ centered at the origin. Then the results in  Melrose and Sj\"ostrand (see  \cite{MS78}, and also \cite{MS82})
imply the following:
\begin{lem} \label{l.sm1}
Let $n\ge2$.
Suppose that $\Omega$ satisfies Assumption {\rm A}.
Let $R > 1$  be fixed such that the obstacle $\mathcal{O}$ is inside $B_{R}$. For $\delta>0$ 
and $\rho \in S^*(\overline{\Omega})$ with $d(\rho, \Sigma) \geq \delta$, 
let $t_{\delta,R,\rho}$ be the sojourn time.
Then  for any $N \geq 1$ there exists a constant $C=C_{N,R,\delta}>0$ such that for any  initial data
 $(u_0,u_1) \in  H^1_0(\Omega) \times L^2(\Omega)$ with
\[
\mathrm{supp}\, u_0 \cup \mathrm{ supp}\, u_1 \subseteq \Omega\, \cap B_{R+1},
\]
\[ \mathrm{WF}_b(u_0) \cup \mathrm{WF}_b(u_1) \subseteq \{\rho \in S^*(\overline{\Omega}):d(\rho, \Sigma) \geq \delta \},
\]
the initial-boundary value problem \eqref{eq.1} admits a unique finite energy solution $u(t,x)$ satisfying
the smoothing estimate
 \begin{equation} \label{EQ:aim}
 \|\nabla u(t)\|_{H^N(\Omega \, \cap B_{R+1})} +\|\partial_t u(t)\|_{H^N(\Omega \, \cap B_{R+1})}
\leq C \left( \|\nabla u_0\|_{L^2(\Omega)} + \| u_1\|_{L^2(\Omega)} \right)
 \end{equation}
for any $t \in(t_{\delta,R,\rho}, 2t_{\delta,R,\rho})$.
\end{lem}
\begin{proof}
Theorem (0.11) in \cite{MS78} states that if $\Omega$ satisfies Assumption {\rm A}, then
$\mathrm{WF}_b(u)$ propagates along the bicharacteristic flow.
 It is well known that
for any cut-off function $\chi$, the operator $$ \chi e^{\im t \sqrt{-\left. \Delta \right|_{D}}}\chi$$ has a smoothing property
(see, e.g., Vodev \cite{Vodev2002}), i.e.,
it maps $ f \in L^2(\Omega) $ into
\[
\chi e^{\im t \sqrt{-\left. \Delta \right|_{D}}}\chi f \in H^{N}(\Omega),
\]
provided that the wave front set of $f$ has no intersection with the trapped set. Hence, we get the
estimate \eqref{EQ:aim}.
The proof of Lemma \ref{l.sm1} is complete.
\end{proof}

Based on Lemma \ref{l.sm1}, it is possible to consider the smoothing estimate for the inhomogeneous problem:
\begin{equation}\label{eq.7}
\left\{
\begin{aligned}
& \partial_{t}^2 u - \left. \Delta\right\vert_D u = F(t,x),  & \quad (t,x)\in \R\times \Omega,\\
   & u(0,x)=u_0(x), \quad \partial_t u(0,x) = u_1(x),& \quad x\in \Omega.
\end{aligned}\right.
\end{equation}
More precisely, by using Duhamel principle, we have the following:

\begin{lem} \label{l.sm2}
Let $n\ge2$. Suppose that $\Omega$ satisfies Assumption {\rm A}.
Let $\delta $ and $R$ be as in Lemma {\rm \ref{l.sm1}}, and let $t_{\delta,R,\rho}$ be the sojourn time associated to the homogeneous problem \eqref{eq.1}
with initial data
 $(u_0,u_1) \in  H^1_0(\Omega) \times L^2(\Omega)$ satisfying
\[
\mathrm{supp}\, u_0 \cup \mathrm{ supp} \, u_1 \subseteq \Omega \cap B_{R+1},
\]
\[ \mathrm{WF}_b(u_0) \cup \mathrm{WF}_b(u_1) \subseteq \{\rho \in S^*(\overline{\Omega}) : d(\rho, \Sigma) \geq \delta \}.
\]
Then  for any $F \in C(\R_+; H^N(\Omega) \cap H^1_0(\Omega))$ with $N\ge1$,
the solution $u$ to the initial-boundary value problem \eqref{eq.7} satisfies
\begin{align*}
 & \| \nabla u(t)\|_{H^N(\Omega\, \cap B_{R+1})}
 +\|\partial_t  u(t)\|_{H^{N}(\Omega\, \cap B_{R+1})} \\
\leq &\,  C  \Big( \|\nabla u_0
\|_{L^2(\Omega)}+ \| u_1\|_{L^2(\Omega)}
+\|F\|_{L^1\big((0,t);H^N(\Omega)\big)} \Big)
\end{align*}
for any $t \in (t_{\delta,R,\rho},2t_{\delta,R,\rho})$.

In addition, if we suppose that
\[
\mathrm{supp}_x\, F(t,\cdot)\subseteq \Omega \cap B_{R+1}
\]
and
\[  \mathrm{WF}_b(F(t,\cdot) ) \subseteq \{\rho \in S^*(\overline{\Omega}) : d(\rho, \Sigma) \geq \delta \}
\]
for any $t \in (t_{\delta,R,\rho},2t_{\delta,R,\rho})$, then 
\begin{equation} \label{psr1}
\begin{split}
 & \| \nabla u(t)\|_{H^N(\Omega\, \cap B_{R+1})}
 +\|\partial_t  u(t)\|_{H^{N}(\Omega\, \cap B_{R+1})} \\
\leq &\,  C  \Big( \|\nabla u_0
\|_{L^2(\Omega)}+ \| u_1\|_{L^2(\Omega)}
+\|F\|_{L^1\big((0,t);L^2(\Omega)\big)} \Big)
\end{split}
\end{equation}
for any $t \in (t_{\delta,R,\rho},2t_{\delta,R,\rho})$.
\end{lem}


\section{Resolvent estimates in high frequency}
\label{sec:4}

Our aim in this section is to obtain a priori estimate for a solution
$w=w_\ep$ to the Helmholtz equation of the form
\begin{equation} \label{eq.2.1}
  \left(  \left. -\Delta \right\vert_D - (\lambda^2 \pm \im \varepsilon) \right)w =f, \quad \lambda>0, \quad \ep>0.
\end{equation}
For this purpose, we need to find several estimates near the boundary and those outside a large ball separately.
We recall the notation
\[
B_R=\{x\in \Rn:|x|<R\},
\]
and choose $R>1$ so large that the obstacle $\mathcal{O}$ is inside the ball $B_R.$ 

We shall prove here the following:
\begin{prop} \label{lem:lem 4.1}
Let $n\ge2$. Suppose that $\Omega$ satisfies Assumptions {\rm A}, {\rm B} and {\rm C}.
Let $s>1/2$. Assume that $\langle x\rangle^s f$ belongs
to $H^{M+\alpha}(\Omega)$ for some $M\geq 0$ and $\alpha>(n-1)/2$. Then there exists a sufficiently small $\ep_0>0$ such that
the solution $w=w_\ep$ to the Helmholtz equation \eqref{eq.2.1} satisfies the following
estimate{\rm :}
\begin{equation} \label{eq.sl1}
\begin{split}
    & \lambda \| w \|_{H^M(\Omega\, \cap B_R)}   +  \| \nabla w \|_{H^M(\Omega\, \cap B_R)} \\
        & \qquad \qquad
        +\lambda \|\langle x \rangle^{-s} w \|_{H^M(|x|>R)}  + \|\langle x \rangle^{-s}\nabla w \|_{H^M(|x|>R)}\\
 \lesssim&\,
(1+\lambda)^{N_0} \|\langle x \rangle^s  f\|_{H^{M+\alpha}(\Omega)}
 \end{split}
\end{equation}
for  $\lambda \gg 1$ and $\ep\in(0,\ep_0]$, where $N_0$ is the constant appearing in Assumption {\rm C}.
\end{prop}

We begin by investigating the relation between the exterior norm and interior one of $w$.
The desired estimate \eqref{eq.sl1} suggests to introduce the following interior and exterior norms:
\begin{equation}  \label{inex1} \left\{
\begin{aligned}
 \vertiii{w}_{\mathrm{int},\lambda,R,M}
&:= \lambda \| w \|_{H^M(\Omega\, \cap B_R)}   +  \| \nabla w \|_{H^M(\Omega\, \cap B_R)},  \\
 \vertiii{w}_{\mathrm{ext},s,\lambda, R,M} &:=  \lambda \|\langle x \rangle^{-s} w \|_{H^M(|x| >R)}
+ \|\langle x \rangle^{-s}\nabla w \|_{H^M(|x| > R)}. 
\end{aligned} \right.
\end{equation}
Then we have the following:
\begin{lem} \label{lem:lem 4.2}
Let $n\ge2$. Assume that $\langle x \rangle^s f\in H^M(\Omega)$ for some $s>1/2$ and $M\geq0$.
Then the solution $w=w_\ep$ to the Helmholtz equation \eqref{eq.2.1} satisfies the following
inequality{\rm :}
\begin{equation} \label{eq.sl1mmaa}
\vertiii{w}_{\mathrm{ext},s,\lambda, R,M}
 \lesssim   \, \|\langle x \rangle^s f \|_{H^M(|x| > R)}
 +\gamma \vertiii{w}_{\mathrm{int}, \lambda,R,M}
\end{equation}
for $\lambda\gg1$ and $\ep \in (0,1)$, where $\gamma>0$ is arbitrarily small.
\end{lem}
\begin{proof} It is sufficient to prove \eqref{eq.sl1mmaa} for $1/2<s\leq 1$. 
We have only to show that
\begin{equation}\label{eq.sl1mm}
\begin{split}
   & \lambda \|\langle x \rangle^{-s} w \|_{H^M(|x| >R)}  + \|\langle x \rangle^{-s}\nabla w \|_{H^M(|x| > R)} \\
 \lesssim  & \, \|\langle x \rangle^s f \|_{H^M(|x| > R)}
 +\gamma \left(\lambda \| w \|_{H^M(R-1 < |x| < R))} + \| \nabla w \|_{H^M(R-1< |x| < R)}\right)
\end{split}
\end{equation}
for $\lambda\gg1$ and $\ep \in (0,1)$, and for an arbitrary small $\gamma>0$.
Indeed, by using Proposition 2.4 from Cardoso and Vodev \cite{VoCa2002}, we deduce that
\begin{equation} \label{eq.sl11}
\begin{split}
   & \lambda \|\langle x \rangle^{-s} w \|_{L^2(|x| >{R})}  + \|\langle x \rangle^{-s}\nabla w \|_{L^2(|x| >{R})} \\
 \lesssim  & \, \|\langle x \rangle^s f \|_{L^2(|x| >{R})} + \frac{1}{\sqrt{\lambda}}  \|\partial_r w \|^{\frac12}_{L^2(|x| =R)}
 \| w \|^{\frac12}_{L^2(|x| =R)}
 \\
&\, +\gamma \left(\lambda \| w \|_{L^2(R-1 < |x| <{R})} +   \| \nabla w \|_{L^2(R-1< |x| <{R})}\right)
\end{split}
\end{equation}
for any $\ep\in(0,1)$ and an arbitrarily small $\gamma >0$.
Resorting to the trace theorem and interpolation inequalities, we estimate $L^2$-norms of $\pa_r w$ and $w$ over the sphere $\{|x|=R\}$ as
\[
\left\{
\begin{gathered}
\|\partial_r w \|_{L^2(|x| =R)} \lesssim \|w \|_{H^{\frac32}( R < |x| < R+1)}
 \le  \|w \|^{\frac12}_{H^{2}(  R < |x| < R+1)}
\|w \|^{\frac12}_{H^{1}(  R < |x| < R+1)},\\
\| w \|_{L^2(|x| = R)} \lesssim \|w \|_{H^{\frac12}(  R < |x| < R+1)}
 \le  \|w \|^{\frac12}_{H^{1}(  R < |x| < R+1)}\|w \|^{\frac12}_{L^2(R < |x| < R+1)},
\end{gathered}
\right.
\]
respectively.
Here we have, by using equation \eqref{eq.2.1},
\[
 \|w \|^{\frac12}_{H^2(  R < |x| < R+1)} \lesssim \lambda \|w \|^{\frac12}_{L^2(  R < |x| < R+1)}
+  \|f \|^{\frac12}_{L^2(  R < |x| < R+1)}.
\]
In this way we deduce from \eqref{eq.sl11} that
\begin{align*}
   & \lambda \|\langle x \rangle^{-s} w \|_{L^2(|x| >R)}  + \|\langle x \rangle^{-s}\nabla w \|_{L^2(|x| >R)} \\
 \lesssim  & \, \|\langle x \rangle^s f \|_{L^2(|x| >R)}
 + \frac{1}{\sqrt{\lambda}} \left( \lambda \| w \|_{L^2( R < |x| < R+1)}  + \|\nabla w \|_{L^2( R < |x| < R+1)} \right)\\
&  +\gamma\left(\lambda \| w \|_{L^2(R-1 < |x| <  R))} +   \| \nabla w \|_{L^2(R-1 <|x| < R)}\right)
\end{align*}
and moreover, by using the fact  that
\begin{equation} \label{eq.classical}
\sum_{|\alpha| \leq M}\| \partial_x^\alpha w \|_{L^2(|x| >R)} \sim \Big\|\left( 1  \left. -\Delta \right\vert_D \right)^{\frac{M}{2}} w \Big\|_{L^2(|x| >R)}
\end{equation}
and
\begin{equation} \label{eq.comm}
\sum_{|\alpha| \leq M}\|\langle x \rangle^{-s} \partial_x^\alpha w \|_{L^2(|x| >R)} \sim \Big\|\langle x \rangle^{-s} \left( 1  \left. -\Delta \right\vert_D \right)^{\frac{M}{2}} w \Big\|_{L^2(|x| >R)},
\end{equation}
and by using the property that $ 1  \left. -\Delta \right\vert_D $ commutes with the operator 
$ \left. -\Delta \right\vert_D - (\lambda^2 \pm \im \varepsilon)$ 
in the equation \eqref{eq.2.1},
we find that
\begin{equation}\label{eq.res}
\begin{split}
   & \lambda \|\langle x \rangle^{-s} w \|_{H^M(|x| >R)}  + \|\langle x \rangle^{-s}\nabla w \|_{H^M(|x| >R)} \\
 \lesssim  & \, \|\langle x \rangle^s f \|_{H^M(|x| >R)}
 + \frac{1}{\sqrt{\lambda}} \left( \lambda \| w \|_{H^M( R < |x| < R+1)}  + \|\nabla w \|_{H^M( R < |x| < R+1)} \right)\\
& +\gamma\left(\lambda \| w \|_{H^M(R-1<|x| <  R))} +   \| \nabla w \|_{H^M(R-1 < |x| < R)}\right).
\end{split}
\end{equation}
Here, \eqref{eq.classical} is a consequence of the elliptic regularity theorem, and \eqref{eq.comm} is proved
in Lemma \ref{lem:A5} of Appendix \ref{App:Appendix A}.
Thus, taking  $\lambda $ sufficiently large in \eqref{eq.res}, we  obtain \eqref{eq.sl1mm}.
The proof of Lemma \ref{lem:lem 4.2} is now complete.
\end{proof}

Our next  step is to control microlocally $H^M$-norm of $w$ close to the trapped set, i.e., in $\delta$-neighbourhood $\Sigma_\delta$ of the 
trapped set $\Sigma$. For this purpose,
we consider a
properly supported pseudo-differential operator $P^\delta_{\mathrm{tr}}=P^\delta_{\mathrm{tr}}(x,D_x)$
of order $0$ with symbol $p^\delta_{\mathrm{tr}}(x,\xi)$ satisfying
\begin{equation} \label{EQ:proper}
\left\{
\begin{aligned}
&  \mathrm{supp}_x \, p^\delta_{\mathrm{tr}}(x,\xi) \subseteq \overline{\Omega} \cap B_{R+2}, \\
&  \mathrm{supp}_{x,\xi} \, p^\delta_{\mathrm{tr}}(x,\xi) \subseteq \left\{ (x,\xi) \in \Sigma_{2\delta} : x \in \overline{\Omega} \cap   B_{R+2}  \right\},\\
&  p^\delta_{\mathrm{tr}}(x,\xi)=1 \quad \text{for $ (x,\xi) \in \Sigma_{\delta}$ \ with $ x \in  \overline{\Omega} \cap B_{R+1}$.}
\end{aligned}
\right.
\end{equation}
Then  the solution $w$ to the Helmholtz equation \eqref{eq.2.1} can be decomposed into
\begin{equation*}
w = w_{\mathrm{I}} + w_{\mathrm{II}},
\end{equation*}
where $w_{\mathrm{I}}$ satisfies the equation
\begin{equation}\label{eq.w11}
\left(\left. -\Delta\right|_D-(\lambda^2\pm \im \ep) \right)w_{\mathrm{I}}=P^\delta_{\mathrm{tr}} f,
\end{equation}
while $w_{\mathrm{II}}$ solves
\begin{equation}\label{eq-wii}
\left(\left. -\Delta\right|_D-(\lambda^2\pm \im \ep)\right)w_{{\mathrm{II}}}=\left(1-P^\delta_{\mathrm{tr}}\right) f.
\end{equation}
Since the support of $p^\delta_{\mathrm{tr}}(x,\xi)$ is compact in $\overline{\Omega} \cap B_{R+2}$,
$P^\delta_{\mathrm{tr}} f$ is compactly supported in $\overline{\Omega} \cap B_{R+2}$.
Hence, applying Theorem 2 from Burq \cite{Burq-Acta} (see also Vodev \cite{Vodev2000}), 
we have
\begin{equation*}
     \lambda \|w_{\mathrm{I}}\|_{L^2(\Omega\, \cap B_R)}
     + \|\nabla w_{\mathrm{I}}\|_{L^2(\Omega\, \cap B_R)}
     \le Ce^{c \lambda}\|P^\delta_{\mathrm{tr}} f\|_{L^2(\Omega)}.
\end{equation*}
This estimate is established for $0<\varepsilon \leq C e^{-c\lambda}$ with $\lambda \ge \lambda_0\, (\gg1)$ in \cite{Burq-Acta}.
However, one can easily extend  it for $ C e^{-c\lambda} \leq \varepsilon \leq \ep_0,$ where
\[
\ep_0=Ce^{-c\lambda_0},
\]
by using the classical estimate
\begin{equation*} 
 \left\|(-\Delta|_{D}-(\lambda^2 \pm \im \ep))^{-1} g\right\|_{L^2(\Omega)} \lesssim \frac{\|g\|_{L^2(\Omega)}}{\ep}.
\end{equation*}
Thus we find that
\begin{equation*} 
     \lambda \|w_{\mathrm{I}}\|_{H^M(\Omega\, \cap B_R)}
     + \|\nabla w_{\mathrm{I}}\|_{H^M(\Omega\, \cap B_R)}
     \lesssim e^{c \lambda}\|P^\delta_{\mathrm{tr}} f\|_{H^M(\Omega)}
\end{equation*}
for any $M\geq0$, $\lambda\ge \lambda_0$ and $\ep \in(0,\ep_0]$.
Therefore, the cut-off resolvent operator
\[
\mathds{1}_{\Omega \, \cap B_ R} \, \left(-\Delta|_D-z^2\right)^{-1} \,  \mathds{1}_{\Omega \, \cap B_ R}
\]
possesses no poles in the domain
\[
\{z\in \mathbb{C}: |\mathrm{Re}(z)|\ge \lambda_0, \, |\mathrm{Im}(z^2)|\le \ep_0\},
\]
where $\mathds{1}_{\Omega \, \cap B_ R}$ stands for the characteristic function of the set $\Omega \, \cap B_ R$.
Hereafter, we fix $\ep_0>0$ and $\lambda_0>1$.
\\

Summarizing the above argument, we have:
\begin{lem} \label{lem:lem 4.3}
Let $n\ge2$. Assume that $f\in H^M(\Omega)$ for some $M\geq 0$.
Then there exists a constant $c>0$ such that the solution $w_{\mathrm{I}}$ to the elliptic equation \eqref{eq.w11} satisfies the following
estimate{\rm :}
\begin{equation} \label{eq2.88amn}
     \lambda \|w_{\mathrm{I}}\|_{H^M(\Omega\, \cap B_R)}
     + \|\nabla w_{\mathrm{I}}\|_{H^M(\Omega\, \cap B_R)} 
     \lesssim e^{c \lambda}\|P^\delta_{\mathrm{tr}} f\|_{H^M(\Omega)}
\end{equation}
for any $\lambda\geq \lambda_0$, $\ep \in(0,\ep_0]$, and for any properly supported pseudo-differential operator $P^\delta_{\mathrm{tr}}(x,D_x)$ of order $0$ with symbol $p^\delta_{\mathrm{tr}}(x,\xi)$ satisfying
\eqref{EQ:proper}.
\end{lem}

We have to estimate $\|P^\delta_{\mathrm{tr}} f\|_{H^M(\Omega)}.$
In fact, we have a sharper result as follows:
\begin{lem} \label{lem:Pf}
Let $n\geq 2$.   Assume that $f \in H^{M+\alpha}(\Omega)$ for some $M\geq 0$ and $\alpha>(n-1)/2$. 
Let $\Sigma_\delta$ be the $\delta$-neighbourhood of the trapped set $\Sigma$ having a covering of type \eqref{eq.cove}, so that Assumption~{\rm B} is satisfied. Then
\begin{equation} \label{EQ:22dl}
\left\| \langle x \rangle^s P^\delta_{\mathrm{tr}}f \right\|_{H^{M}(\Omega)}  \lesssim
 \delta^{\frac{n-1-\nu}{2}}
\|  f\|_{H^{M+\alpha}(\Omega)} 
\qquad 
\end{equation}
for some $s>1/2$ and $\nu\in[0,n-1)$. 
\end{lem}
\begin{proof}
Noting from the definition of the properly supported pseudo-differential operator $P^\delta_{\mathrm{tr}}(x,D)$, 
we have only to show that 
\begin{equation} \label{EQ:2dl}
\left\| P^\delta_{\mathrm{tr}}f \right\|_{H^{M}(\Omega)}  \lesssim
 \delta^{\frac{n-1-\nu}{2}}
\|  f\|_{H^{M+\alpha}(\Omega \, \cap  B_{R_0})}
\end{equation}
for $\nu\in[0,n-1)$, where $R_0>0$ is such that the ball of radius $R_0$ covers the union of compact sets 
$K_{\delta,j}$, $j = 1,\ldots, N(\delta)$, in Assumption B. 
Moreover, it is sufficient to prove \eqref{EQ:2dl} for $M=0$, since the case $M>0$
is similar to that of $M=0$.

Our Assumption B guarantees that we can find a real $\delta_0>0$ such that for $\delta < \delta_0$,  
a covering of $\Sigma_\delta$ is taken as
$$ \Sigma_\delta = \bigcup_{j=1}^{N(\delta)} \left(K_{\delta,j} \times E_{\delta,j} \right),$$
where $N(\delta)$ is the number satisfying
\begin{equation*}
  N(\delta) \sim \delta^{-\nu}, \quad 0 \leq \nu <n-1.
\end{equation*}
We consider the operators 
$P^{\delta}_{\mathrm{tr},j}(x,D_x)$ for 
$j=1,\ldots,N(\delta)$,
which are restrictions of $P^{\delta}_{\mathrm{tr}}(x,D_x)$ to $K_{\delta,j} \times E_{\delta,j}$.
By using the standard extension operator
$$
\begin{array}{ccc}
H^M(\Omega) & \longrightarrow & H^M(\mathbb{R}^n) \\
\rotatebox{90}{$\in$} & & \rotatebox{90}{$\in$} \\
f & \longmapsto &\widetilde{f},
\end{array}
$$
we can write
$$  P^{\delta}_{\mathrm{tr},j}(x,D_x)f  = \left.  \widetilde{P}^{\delta}_{\mathrm{tr},j}(x,D_x) \widetilde{f}  \right|_{\Omega},$$
where we put 
\begin{equation} \label{eq.tr}
\widetilde{P}^{\delta}_{\mathrm{tr},j}(x,D_x) := \Phi(x) \psi^{(j)}_\delta(D_x) \varphi^{(j)}_{r(\delta)} (x).
\end{equation}
Here, $\Phi(x)$ is compactly supported smooth function on $\Rn$, and is independent of $j$. 
Each $\va^{(j)}_{r(\delta)}(x)$ is a smooth function on
$\R^n$ and its support has the diameter $r(\delta)$, while each $\psi^{(j)}_\delta(D_x)$ is the pseudo-differential operator of order $0$ whose symbol
$\psi^{(j)}_\delta(\xi)$ is supported in $E_{2\delta,j}$ and equal to 1 on $E_{\delta,j}.$
We recall that each $E_{2\delta,j}$ is a ball  on the sphere $\mathbb{S}^{n-1}$  of radius $2\delta$ near the center of $E_{\delta,j}$. 
Put $g= \varphi^{(j)}_{r(\delta)} \widetilde{f}$. 
Since $\Phi(x) $ is bounded, we estimate   
\begin{align*}
\big\|\Phi(\cdot) \psi^{(j)}_\delta(D_x) g \big\|_{L^2(\Rn)}\lesssim &\, \big\| \psi^{(j)}_\delta(D_x) g \big\|_{L^2(\Rn)}\\
 \lesssim &\, \mu(E_{\delta,j})^{\frac12} \left(\int_0^\infty\|\widehat{g}(\rho \, \cdot) 
\|^2_{L^\infty(\mathbb{S}^{n-1})} \rho^{n-1} d\rho\right)^{\frac12},
\end{align*} 
where $\mu(E_{\delta,j})$ stands for the $(n-1)$-dimensional Lebesgue measure of 
$E_{\delta,j}$, and hence, $\mu(E_{\delta,j})\sim \delta^{n-1}$. 
By the Sobolev embedding theorem on $\mathbb{S}^{n-1}$, we see from Plancherel's identity that 
\begin{align*}
\int^\infty_0 \|\widehat{g}(\rho \, \cdot) \|^2_{L^\infty(\mathbb{S}^{n-1})} \rho^{n-1}\, d\rho
\lesssim&\, \int^\infty_0  \|\widehat{g}(\rho \, \cdot) \|^2_{H^\alpha(\mathbb{S}^{n-1})} \rho^{n-1}\, d\rho\\
=&\, \int^\infty_0 \|g(r \, \cdot) \|^2_{H^\alpha(\mathbb{S}^{n-1})}r^{n-1}\, dr
\end{align*}
for $\alpha>(n-1)/2$. Hence, by using these estimates with $g= \varphi^{(j)}_{r(\delta)}(\cdot) \widetilde{f}$, we get
\begin{align*}
\big\|\psi^{(j)}_\delta(D_x) \varphi^{(j)}_{r(\delta)}(\cdot) \widetilde{f} \big\|_{L^2(\Rn)} \lesssim&\, \delta^{\frac{n-1}{2}} 
\left(\int_0^\infty \big\|  \varphi^{(j)}_{r(\delta)}(r\cdot) f(r \cdot) \big\|^2_{H^\alpha(\mathbb{S}^{n-1})} r^{n-1} dr \right)^{\frac12}\\
\lesssim &\, \delta^{\frac{n-1}{2}} \big\|\varphi^{(j)}_{r(\delta)}f \big\|_{H^\alpha(\Omega)}.
\end{align*}
Since each $\va^{(j)}_{r(\delta)}$ is supported in the ball with radius $R_0$ by Assumption B, we use  this estimate with \eqref{eq.tr}, combined with
\begin{equation} \label{eq.dis}
\Big\|\sum_j h_j \Big\|^2_{L^2(\mathbb{S}^{n-1})} \sim \sum_j \| h_j \|^2_{L^2(\mathbb{S}^{n-1})}
\end{equation}
for $h_j$ with disjoint (almost disjoint ) supports in $\mathbb{S}^{n-1}$. In this way, taking sum over $j=1,\ldots,N(\delta)$, we conclude that 
$$ \big\|  \widetilde{P}^\delta_{\mathrm{tr}}\widetilde{f} \big\|_{L^2(\R^n)} \lesssim \delta^{\frac{n-1}{2}}\cdot
 \delta^{-\frac{\nu}{2}}
\left\| f \right\|_{H^\alpha(\Omega \, \cap B_{R_0})} .
\qquad 
$$ 
Indeed, if the supports of non-negative functions $h_j$ are almost disjoint (which means that there is a fixed integer $k$ so that for any $\omega\in \mathbb{S}^{n-1}$ there are at most $k$ functions $h_j$ such that $h_j(\omega)\neq 0$), then we have for any $p\in(1,\infty)$, 
$$ \Big|\sum_j  h_j(\omega) \Big|^p \sim \sum |h_j(\omega)|^p,$$
which implies \eqref{eq.dis}.
Thus we get 
the required estimate
\eqref{EQ:2dl} for $M=0$. The proof of Lemma \ref{lem:Pf} is complete.
\end{proof}

As a consequence of Lemmas \ref{lem:lem 4.3} and \ref{lem:Pf}, we have:
\begin{lem} \label{lem:lem 4.4}
Let $n\ge2$.
Suppose that $\Omega$ satisfies Assumption {\rm B}.
For any $\lambda\geq \lambda_0$ and $\kappa \in (0,1) $, take $\delta \in (0,1)$ so that
\begin{equation}\label{dlr1}
     \frac{n-1-\nu}{2}\log \left(\frac{1}{\delta}\right)  \gtrsim \lambda + \log\left(\frac{1}{\kappa}\right).
\end{equation}
Then the solution $w_{\mathrm{I}}=w_{\mathrm{I}}^{(\delta)}$ to \eqref{eq.w11} satisfies the estimate
\begin{equation}\label{Vbest}
    \lambda \|w_{\mathrm{I}}\|_{H^M(\Omega\, \cap B_R)}
     + \|\nabla w_{\mathrm{I}}\|_{H^M(\Omega\, \cap B_R)}
     \leq \kappa \|  f\|_{H^{M+\alpha}(\Omega)}
\end{equation}
for any  $f \in H^{M+\alpha}(\Omega)$ with $M\geq0$ and $\alpha>(n-1)/2$. 
\end{lem}
\begin{proof}
Applying \eqref{EQ:22dl} in Lemma \ref{lem:Pf} to \eqref{eq2.88amn} in Lemma \ref{lem:lem 4.3}, we get
\[
   \lambda \|w_{\mathrm{I}}\|_{H^M(\Omega\, \cap B_R)}
     + \|\nabla w_{\mathrm{I}}\|_{H^M(\Omega\, \cap B_R)}
     \lesssim e^{c\lambda} \delta^{\frac{n-1-\nu}{2}}
     \| f\|_{H^{M+\alpha}(\Omega)}.
\]
Hence we have only to show that for any fixed $\kappa \in (0,1)$ and $\lambda\geq \lambda_0$,
\begin{equation*} 
e^{c\lambda} \delta^{\frac{n-1-\nu}{2}} \leq \kappa,
\end{equation*}
which is possible due to \eqref{dlr1}. Thus we conclude the estimate \eqref{Vbest}.
The proof of Lemma~\ref{lem:lem 4.4} is complete.
\end{proof}

\begin{remark}
It is crucial to clarify the relation between $\delta$ and $\lambda.$ The inequality \eqref{dlr1} suggests to take equality there and hence
\begin{equation}\label{rdl1}
     \lambda \sim \log \left(\frac{1}{\delta}\right).
\end{equation}
\end{remark}

We are now in a position to derive the estimate for $w_I$. Recalling the notation of the interior norm and exterior 
one (see \eqref{inex1}), we have the following: 
\begin{lem} \label{lem:lem4-6}
Let $n\ge2$. Suppose that $\Omega$ satisfies Assumption {\rm B}.
Assume that $f \in H^{M+\alpha}(\Omega)$ for some $M\geq0$ and $\alpha>(n-1)/2$.
Then the solution $w_{\mathrm{I}}$ to \eqref{eq.w11} satisfies the following estimate{\rm :}
\begin{equation}\label{eq.inter}
\vertiii{w_{\mathrm{I}}}_{\mathrm{int},\lambda, R, M} +\vertiii{w_{\mathrm{I}}}_{\mathrm{ext},s,\lambda, R, M}
\lesssim \|f\|_{H^{M+\alpha}(\Omega)}
\end{equation}
for any $\lambda\geq \lambda_0$, $\ep\in(0,\ep_0]$ and $s>1/2$.
\end{lem}
\begin{proof}
By using the inequality \eqref{eq.sl1mmaa} from Lemma \ref{lem:lem 4.2} and \eqref{EQ:22dl},
we find that
\begin{align*}
\vertiii{w_{\mathrm{I}}}_{\mathrm{ext},s,\lambda, R, M}
 \lesssim &\,  \|\langle x \rangle^s P^{\delta}_{\mathrm{tr}} f \|_{H^M(|x| > R)}
 +\gamma \vertiii{w_{\mathrm{I}}}_{\mathrm{int},\lambda, R, M}\\
 \lesssim &\, \|f \|_{H^{M+\alpha}(\Omega)}
 +\gamma \vertiii{w_{\mathrm{I}}}_{\mathrm{int},\lambda, R, M}
\end{align*}
for $M\geq0$, $\alpha>(n-1)/2$, $\lambda\geq \lambda_0$, $\ep\in(0,\ep_0]$ and $s>1/2$, and for an arbitrarily small $\gamma>0$. On the other hand, it follows from
\eqref{Vbest} in Lemma \ref{lem:lem 4.4} that
\[
\vertiii{w_{\mathrm{I}}}_{\mathrm{int},\lambda, R, M}
 \lesssim \kappa \|f \|_{H^{M+\alpha}(\Omega)}
\]
for $M\geq 0$, $\alpha>(n-1)/2$, $\lambda\geq \lambda_0$, $\ep\in(0,\ep_0]$,
and for any $\kappa\in(0,1)$.
Hence the estimate \eqref{eq.inter} is an immediate consequence of these inequalities. This ends the proof of
Lemma \ref{lem:lem4-6}.
\end{proof}

Let us turn to estimate the solution $w_{{\mathrm{II}}}$  to \eqref{eq-wii}. 

\begin{lem} \label{lem:lem 4.7}
Let $n\ge2$ and $M\geq1$.
Suppose that $\Omega$ satisfies Assumptions {\rm A}, {\rm B} and {\rm C}.
Assume that $\langle x \rangle^s f \in H^{1}(\Omega)$  for some  $s> 1/2$.
For $\lambda>\lambda_0$, take $\delta>0$ so that \eqref{rdl1} holds. Then 
the solution $w_{{\mathrm{II}}}$  to \eqref{eq-wii} satisfies the following estimate{\rm :}
\begin{equation}\label{IIVbest}
\begin{split}
   &  \vertiii{ w_{{\mathrm{II}}}}_{\mathrm{int},\lambda, R,M} \lesssim (1+\lambda)^{N_0}  \|\langle x \rangle^s f \|_{H^1(\Omega)},
\end{split}
\end{equation}
where $N_0$ is the constant appearing in Assumption {\rm C}.
\end{lem}
\begin{proof}
Since \eqref{eq-wii} is an elliptic equation, it follows from the definition of $P^{\delta}_{\mathrm{tr}}$ that 
\begin{equation} \label{eq.wII1}
\mathrm{WF}_b(w_{{\mathrm{II}}}) \subseteq \mathrm{WF}_b((1-P^{\delta}_{\mathrm{tr}})f) \subseteq
\{\rho \in S^*(\overline{\Omega}) : d(\rho, \Sigma) \geq \delta \}. 
\quad
\end{equation}
Let us take a spectral parameter $z = z(\lambda, \varepsilon)$
such that
\[
\text{$ z^2 = \lambda^2 \pm \im \varepsilon$,}
\]
where $z$ is defined by taking the principal branch of logarithm, i.e.,
\[
z = e^{\frac12 \mathrm{Log}(\lambda^2 \pm \im \varepsilon)}.
\]
For $\lambda \gg 1$ and $\varepsilon >0$ small enough, we have
\begin{equation*} 
z = \sqrt{\lambda^2 \pm \im \varepsilon} =
  \lambda \pm \frac{\im \varepsilon}{2\lambda}
+ O\left(\frac{\ep^2}{\lambda^3}\right).
\end{equation*}
Furthermore, taking a non-negative smooth function $\chi(r)$, $r=|x|$, such that
\[
 \chi(r)=\left\{
\begin{aligned}
1, \quad & 0 \le r \le R; \\
0, \quad & r \geq R+1,
\end{aligned}
\right.
\]
we set
\[
U_\pm (t,x) = e^{\pm \im z t} \chi w_{{\mathrm{II}}}.
\]
Since $w_{{\mathrm{II}}}$ satisfies the equation \eqref{eq-wii}, it follows that
$U_\pm$ are solutions to the following initial-boundary value problems
\begin{equation*}
\left\{
\begin{aligned}
 & \partial_{t}^2 U_\pm - \left. \Delta\right\vert_D U_\pm = F_\pm(t,x),  & \quad (t,x)\in (0,\infty)\times \Omega,\\
 & U_\pm(0)=\chi w_{{\mathrm{II}}}, \quad \partial_t U_\pm(0) = \pm \im z \chi  w_{{\mathrm{II}}}, & \quad x\in \Omega,
\end{aligned}\right.
\end{equation*}
respectively,
where
\[
F_\pm(t,x) =  e^{\pm \im z t}\chi  (1-P^{\delta}_{\mathrm{tr}}) f - e^{\pm\im z t} [\left. -\Delta\right|_D, \chi] w_{{\mathrm{II}}}.
\]
We have to estimate $F_\pm(t,x) $. Since
\begin{equation*}
\mathrm{WF}_b(U_\pm(0)) \cup \mathrm{WF}_b(\pa_t U_\pm(0))=
\mathrm{WF}_b( \chi w_{{\mathrm{II}}}),
\end{equation*}
we use the relation \eqref{eq.wII1} to find that  
\[
\mathrm{WF}_b(U_\pm(0)) \cup \mathrm{WF}_b(\pa_t U_\pm(0))
\subseteq  \{\rho \in S^*(\overline{\Omega}) : d(\rho, \Sigma) \geq \delta \},
\]
and
\[
 \mathrm{WF}_b(F_\pm(t,\cdot) ) \subseteq \{\rho \in S^*(\overline{\Omega}) : d(\rho, \Sigma) \geq \delta \}
 \]
for any $t \in \mathbb{R}$.

Now we are ready to apply \eqref{psr1} from Lemma \ref{l.sm2} for $t \in(t_{\delta},2t_{\delta})$, and we get
\begin{equation}
\begin{split}
& \lambda \|w_{{\mathrm{II}}}\|_{H^{M}(\Omega\, \cap B_{R})} +  \|\nabla w_{{\mathrm{II}}}\|_{H^M(\Omega\, \cap B_{R})} \\
    \lesssim & \, \lambda \|w _{{\mathrm{II}}}\|_{L^2(\Omega \,\cap \, (B_{R+1 }\setminus B_R))}
    +\|\nabla w_{{\mathrm{II}}}\|_{L^2(\Omega\, \cap \, (B_{R+1 }\setminus B_R)} \label{EQ:LOC}\\ 
    &+ \|(1-P^{\delta}_{\mathrm{tr}}) f\|_{L^1\left((0,t);L^2(\Omega\, \cap \, B_{R+1})\right)} \\
 &+\|\nabla w_{{\mathrm{II}}}\|_{L^1\left((0,t);L^2(\Omega\, \cap \, (B_{R+1 }\setminus B_R) )\right)}
 + \| w_{{\mathrm{II}}}\|_{L^1\left((0,t);L^2(\Omega\, \cap \, (B_{R+1 }\setminus B_R)) \right)}
  \end{split}
\end{equation}
for any $M\geq 0$ and $t\in(t_{\delta,R,\rho},2t_{\delta,R,\rho})$. We split the norms over $\Omega \cap B_{R+1}$ in \eqref{EQ:LOC} into those on
$\Omega \cap B_R$ and $\Omega \cap \{R<|x|<R+1\}$, respectively.
In this way, we have, by using $L^2(\Omega)$-boundedness of $1-P^{\delta}_{\mathrm{tr}}$,
$$ \vertiii{w_{{\mathrm{II}}}}_{\mathrm{int},\lambda, R,M} \lesssim 
t_{\delta,R,\rho}\vertiii{\widetilde{\chi} w_{{\mathrm{II}}}}_{\mathrm{ext},s,\lambda, R,1} 
+ t_{\delta,R,\rho} \| f \|_{L^2(\Omega)},$$
where $\widetilde{\chi}= \widetilde{\chi}_R$
is equal to 1 in the ball of radius $R+1$, and 0 outside the ball of radius $R+2.$ 
On the other hand, applying the inequality \eqref{eq.sl1mmaa} from Lemma \ref{lem:lem 4.2},
we can write
$$ 
\vertiii{w_{{\mathrm{II}}}}_{\mathrm{ext},s,\lambda, R,1}
 \lesssim   \, \|\langle x \rangle^s f \|_{H^1(|x| > R)}
 +\gamma \vertiii{w_{{\mathrm{II}}}}_{\mathrm{int}, \lambda,R,1}. \quad 
 $$
Hence, we have
\begin{align*}
\vertiii{w_{{\mathrm{II}}}}_{\mathrm{int},\lambda, R,M} 
\lesssim &\, t_{\delta,R,\rho} \left( \|\langle x \rangle^s f \|_{H^1(|x| > R)}
+\gamma \vertiii{w_{{\mathrm{II}}}}_{\mathrm{int}, \lambda,R,1} \right) 
+ t_{\delta,R,\rho} \| f \|_{H^1(\Omega)} \\
\leq &\, 
 t_{\delta,R,\rho} \gamma \vertiii{w_{{\mathrm{II}}}}_{\mathrm{int}, \lambda,R,1} 
  +2 t_{\delta,R,\rho} \|\langle x \rangle^s f \|_{H^1(\Omega)}
 \end{align*}
for an arbitrary small $\gamma>0.$
In this way there exists constants $C_1,C_2>0$ such that 
\begin{equation*} 
\vertiii{w_{{\mathrm{II}}}}_{\mathrm{int},\lambda, R,M} 
\leq C_1\gamma t_{\delta,R,\rho} \vertiii{w_{{\mathrm{II}}}}_{\mathrm{int}, \lambda,R,1}
+ C_2 t_{\delta,R,\rho}  \|\langle x \rangle^s f \|_{H^1(\Omega)}.
\end{equation*}
At this stage, we need in the sequel to unify the estimates of $w_{{\mathrm{I}}}$ and $w_{{\mathrm{II}}}$, and for this we choose equality in
\eqref{dlr1}, so that $\lambda$ is proportional to  $\log(1/\delta).$

In order to have polynomial loss in $\lambda$ we use Assumption C. Namely, 
we have to choose
$$ t_{\delta,R,\rho} \lesssim \mbox{polynomial in $\lambda$} \sim  \left[\log \left( \frac{1}{\delta}\right)\right]^{N_0} \quad \text{and} \quad 
C_1 \gamma t_{\delta,R,\rho} \leq \frac{1}{2},$$ 
and we conclude from \eqref{rdl1} that 
\begin{equation*} 
 \vertiii{w_{{\mathrm{II}}}}_{\mathrm{int},\lambda, R,M} \lesssim 
 (1+\lambda)^{N_0} \|\langle x \rangle^s f \|_{H^1(\Omega)}.
\end{equation*}
The proof of Lemma \ref{lem:lem 4.7} is now complete.
\end{proof}

We are now in a position to estimate $w_{{\mathrm{II}}}$.
\begin{lem} \label{lem:lem II}
Let $n\ge2$, $M\geq1$, $\alpha>(n-1)/2$ and $s>1/2$.
Suppose that $\Omega$ satisfies Assumptions {\rm A}, {\rm B} and {\rm C}. Then
\begin{equation} \label{eq.123}
\vertiii{w_{{\mathrm{II}}}}_{\mathrm{int},\lambda, R, M} + \vertiii{w_{{\mathrm{II}}}}_{\mathrm{ext},s,\lambda, R, M}
         \lesssim (1+\lambda)^{N_0}\|\langle x \rangle^sf\|_{H^{M+\alpha}(\Omega)}
\end{equation}
for any $\lambda\geq \lambda_0$ and $\ep\in(0,\ep_0]$.
\end{lem}
\begin{proof} The estimate \eqref{eq.123}
follows from combination between the estimate
\eqref{eq.sl1mmaa} from Lemma \ref{lem:lem 4.2} and 
the estimate \eqref{IIVbest} for $\vertiii{w_{{\mathrm{II}}}}_{\mathrm{int},\lambda, R, M}$ from Lemma \ref{lem:lem 4.7}.
\end{proof}

\begin{proof}[End of the proof of Proposition {\rm \ref{lem:lem 4.1}}]
Since $w=w_{{\mathrm{I}}} + w_{{\mathrm{II}}}$, we conclude from \eqref{eq.inter} and \eqref{eq.123} that
\begin{equation*}
\vertiii{w}_{\mathrm{int},\lambda, R, M} + \vertiii{w}_{\mathrm{ext},s,\lambda, R, M}
         \lesssim (1+\lambda)^{N_0} \|\langle x \rangle^s f\|_{H^{M+\alpha}(\Omega)}
\end{equation*}
for $M\geq1$, $\alpha>(n-1)/2$, $\lambda\geq \lambda_0$, $\ep \in(0,\ep_0]$ and $s>1/2$.
This proves the estimate \eqref{eq.sl1} for $M\geq1$. 
Finally, the estimate \eqref{eq.sl1} for $M=0$ follows from the equations \eqref{eq.2.1} and 
\eqref{eq.sl1} for $M\geq1$.  
The proof of Proposition \ref{lem:lem 4.1} is now finished.
\end{proof}


\section{Proofs of Theorems \ref{thm:main} and \ref{thm:mainB}}
\label{sec:5}

In this section we prove Theorems \ref{thm:main} and \ref{thm:mainB}.
The following result is concerned with the resolvent estimate in low frequency.
\begin{thm} \label{thm:asymptotic expansion12}
Let $n$ be an odd integer with $n\ge3$, and let $s>1/2$.  Then there exists an
$\varepsilon_0>0$ such that
\begin{enumerate}
\item[i)] the operators
\[
z^{1-k}\langle x\rangle^{-s} \nabla^k (-\Delta|_D-z^2)^{-1}\langle x\rangle^{-s}, \quad k=0,1,
\]
are $\mathcal{B}(L^2(\Omega))$-valued continuous functions on $\{ z \in \mathbb {C}: |z| \leq \varepsilon_0 \}${\rm ;}
\item[ii)] we have the estimate
\begin{equation*} 
\sum_{k=0}^1 |z|^{1-k}\left\|\langle x\rangle^{-s} \nabla^k (-\Delta|_D-z^2)^{-1}\langle x\rangle^{-s}\right\|_{\mathcal{B}(L^2(\Omega))}
\le C
\end{equation*}
on $\{ z \in \mathbb {C}: |z| \leq \varepsilon_0 \}${\rm ;}
\item[iii)] for any cut-off function $\chi$ which is identically $1$ in a neighbourhood of the obstacle $\mathcal{O}$, the operator
\[
\chi(-\Delta|_D-z^2)^{-1}\chi
\]
can be extended as a $\mathcal{B}(L^2(\Omega))$-valued meromorphic function having no poles in
the disk $\{ z \in \mathbb {C}: |z| \leq \varepsilon_0 \}$.

\end{enumerate}
\end{thm}

\begin{remark}
We use the fact that $n \geq 3$ is odd only in the assertion iii), since the free  cut-off resolvent
\[
\chi (-\Delta-z^2)^{-1} \chi
 \]
is not analytic near the origin for $n$ even.
\end{remark}

The proof of Theorem \ref{thm:asymptotic expansion12} 
is rather long, and postponed
until Section~\ref{sec:7} (see also \cite{Georg-Matsuyama}).\\

The following proposition provides the uniform resolvent estimate in intermediate frequency.  It is the well-known property based on the limiting absorption principle and the proof can be found in Mochizuki \cite{Mochizuki-2} (see also \cite{Mochizuki-1}, Kerler \cite{Kerler} and Wilcox~\cite{Wilcox}).
\begin{prop}\label{prop:prrop 63}
Let $n\ge2$, and let $s>1/2$ and $0<a<b$. Then there exists an $\ep_0>0$ such that
\begin{enumerate}
\item[i)] the operators
\[
\langle x\rangle^{-s}\nabla^k (-\Delta|_D-z^2)^{-1}\langle x\rangle^{-s}, \quad k=0,1,
\]
are $\mathcal{B}(L^2(\Omega))$-valued continuous functions on $\Gamma_{a,b}(\ep_0)$, where we put
\[
\Gamma_{a,b}(\ep_0)=\{ z \in \mathbb {C}: a\le |\mathrm{Re}(z)|\le b, \, |\mathrm{Im} (z) | \leq \varepsilon_0 \} \mathrm{;}
\]
\item[ii)] we have the estimate
\begin{equation*}
\sum_{k=0}^1 \left\|\langle x\rangle^{-s}\nabla^k(-\Delta|_D-z^2)^{-1}\langle x\rangle^{-s}\right\|_{\mathcal{B}(L^2(\Omega))}
\le C \quad \text{on $\Gamma_{a,b}(\ep_0)$\rm{;}}
\end{equation*}

\item[iii)] for any cut-off function $\chi$ which is identically $1$ in a neighbourhood of the obstacle $\mathcal{O}$, the operator
\[
\chi (-\Delta|_D-z^2)^{-1}\chi
\]
can be extended as a $\mathcal{B}(L^2(\Omega))$-valued meromorphic function having no poles in $\Gamma_{a,b}(\ep_0)$.

\end{enumerate}
\end{prop}

The following lemma is useful to prove the theorems. 
\begin{lem} \label{lem:lem 5.4}
Let $\ep_0$ be as in Proposition {\rm \ref{lem:lem 4.1}} and let 
$z \in \mathbb {C}$ be such that  
$$z^2=\lambda^2 \pm \im \ep, \quad \lambda  \gg 1, \quad 0 <\ep  \leq \ep_0 .$$
If $f \in \mathrm{Dom}((1 \left. -\Delta \right\vert_D)^k)$ for some integer $k\geq 0$ has a compact support,
then 
\begin{equation}\label{eq.1-m}
\big\|\chi \left(  \left. -\Delta \right\vert_D  - z^2 \right)^{-1} f \big\|_{L^2(\Omega)}
\lesssim (1+\lambda)^{N_0} \big\|  (\left. 1-\Delta \right\vert_D)^k f \big\|_{L^2(\Omega)}
\end{equation}
is equivalent to 
\begin{equation}\label{eq.2-m}
\big\|\chi \left(  \left. -\Delta \right\vert_D  - z^2 \right)^{-1} h \big\|_{L^2(\Omega)}
\lesssim (1+\lambda)^{N_0+2k} \|  h \|_{L^2(\Omega)}
\end{equation}
for any compactly supported function $h\in L^2(\Omega)$. 
\end{lem}
\begin{proof}
Note that \eqref{eq.1-m} is equivalent to 
\begin{equation}\label{eq.3-m}
\big\|\chi \left(  \left. -\Delta \right\vert_D  - z^2 \right)^{-1} (\left. 1-\Delta \right\vert_D)^{-k}h \big\|_{L^2(\Omega)}
\lesssim (1+\lambda)^{N_0} \|  h\|_{L^2(\Omega)}
\end{equation}
for any compactly supported function $h\in L^2(\Omega)$. Hence, it is sufficient to show that \eqref{eq.2-m} 
is equivalent to \eqref{eq.3-m}. It is clear that if $k=0$, then \eqref{eq.2-m} is equivalent to \eqref{eq.3-m}.
We claim that if $k=1$, then 
\eqref{eq.2-m} is equivalent to \eqref{eq.3-m}.
Suppose that \eqref{eq.3-m} holds for $k=1$. 
By using the relation:
\begin{equation}\label{eq.rel}
(z^2+ \left. \Delta \right\vert_D )^{-1} = (1+z^2)(z^2+\left. \Delta \right\vert_D)^{-1}(1 \left. -\Delta \right\vert_D )^{-1}
 -  (1 \left. -\Delta \right\vert_D )^{-1} ,
 \end{equation}
we find that  
\begin{align*}
\big\|\chi \left(  \left. -\Delta \right\vert_D  - z^2 \right)^{-1} h \big\|_{L^2(\Omega)}
 \leq &\, \big\|\chi (1+z^2)(\left. -\Delta \right\vert_D-z^2 )^{-1}(1 \left. -\Delta \right\vert_D )^{-1} h \big\|_{L^2(\Omega)}\\
& 
+\big\|\chi \left(  1\left. -\Delta \right\vert_D \right)^{-1} h \big\|_{L^2(\Omega)}\\
\lesssim &\,  (1+\lambda)^2(1+\lambda)^{N_0} \|h \|_{L^2(\Omega)}
+\|h \|_{L^2(\Omega)}\\
\lesssim &\, (1+\lambda)^{N_0+2} \|  h \|_{L^2(\Omega)},
\end{align*}
which proves \eqref{eq.2-m} for $\kappa=1$. 
Conversely, we suppose that \eqref{eq.2-m} holds for $\kappa=1$. Again by using \eqref{eq.rel}, we have 
\begin{align*}
& (1+\lambda^2) \big\|\chi \left(  \left. -\Delta \right\vert_D  - z^2 \right)^{-1} (\left. 1-\Delta \right\vert_D)^{-1}h \big\|_{L^2(\Omega)}\\
\lesssim &\, \big\|\chi \left(  \left. -\Delta \right\vert_D  - z^2 \right)^{-1} h\big\|_{L^2(\Omega)}
+\big\| (\left. 1-\Delta \right\vert_D)^{-1}h \big\|_{L^2(\Omega)}\\
\lesssim &\, (1+\lambda)^{N_0+2}\| h \|_{L^2(\Omega)}.
\end{align*}
This proves \eqref{eq.3-m} for $\kappa=1$. Consequently, we conclude that 
if $k= 1$, then \eqref{eq.2-m} is equivalent to \eqref{eq.3-m}.

Furthermore, we proceed inductively with respect to an integer $k$. We use the identity
\begin{equation}\label{eq.relk1}
(z^2+ \left. \Delta \right\vert_D )^{-1} (1 \left. -\Delta \right\vert_D )^{-\kappa} = (1+z^2)(z^2+\left. \Delta \right\vert_D)^{-1}(1 \left. -\Delta \right\vert_D )^{-1-k}
 -  (1 \left. -\Delta \right\vert_D )^{-1-k} ,
 \end{equation}
so that using the equivalence between \eqref{eq.1-m} and \eqref{eq.2-m} for $\kappa$ we deduce the equivalence \eqref{eq.1-m} and \eqref{eq.2-m} for $k+1.$ Finally, by interpolation, we conclude that the assertion is valid for 
any integer $k\geq 0$. 
The proof of Lemma \ref{lem:lem 5.4} is complete. 
\end{proof}

Theorem  \ref{thm:main} is a consequence of Proposition  \ref{lem:lem 4.1} and the argument used in 
Tang and Zworski \cite{TZ00}, and Burq \cite{Bu2004}. The detailed proof is as follows.

\begin{proof}[Proof of Theorem {\rm \ref{thm:main}}]

If $\lambda >0$ is bounded, then we can use the estimates of Theorem~\ref{thm:asymptotic expansion12} and Proposition \ref{prop:prrop 63}. For this reason, let us assume that $\lambda \gg 1.$
In order to deduce the weighted $L^2$-estimate \eqref{eq.ME1}, we recall the estimate \eqref{eq.sl1} for $M=0$ from 
Proposition \ref{lem:lem 4.1}:
\begin{equation}\label{eq.C1}
    \lambda \left\| \langle x \rangle^{-s} w \right\|_{L^2(\Omega)} + \left\| \langle x \rangle^{-s} \nabla w \right\|_{L^2(\Omega)} \lesssim (1+\lambda)^{N_0} \|\langle x \rangle^s f\|_{H^{\alpha}(\Omega)}
\end{equation}
for $ s >1/2,$ $\alpha>(n-1)/2$, $\lambda \gg 1$ and $\langle x \rangle^s f \in H^{\alpha}(\Omega)$.
$H^\alpha$-$L^2$-estimate for the weighted resolvent in \eqref{eq.C1} still grows polynomially.  Now, 
by using \eqref{eq.C1} for $\alpha$ restricted to even integer, we write
\begin{equation} \label{eq.gAi1}
   \lambda \big\| \chi \left(  \left. -\Delta \right\vert_D  - z^2 \right)^{-1} f \big\|_{L^2(\Omega)} 
    \lesssim (1+\lambda)^{N_0} \left\|(1 \left. -\Delta \right\vert_D)^{\frac{\alpha}{2}} f \right\|_{L^2}
\end{equation}
for $\lambda\gg 1$ and for any compactly supported function $ f  \in \mathrm{Dom}((1 \left. -\Delta \right\vert_D)^{\alpha/2})$, 
where we put 
\[
z^2=\lambda^2 \pm \im \ep, \quad \lambda \gg 1.
\] 
Then we apply Lemma \ref{lem:lem 5.4} to deduce from \eqref{eq.gAi1} that  
\begin{equation} \label{eq.gAi2}
  \lambda  \big\| \chi \left(  \left. -\Delta \right\vert_D  - z^2 \right)^{-1} f \big\|_{L^2(\Omega)} 
    \lesssim (1+\lambda)^{N_0+\alpha}\|f\|_{L^2(\Omega)}
\end{equation}
for any compactly supported functions $f \in L^2(\Omega).$
However, once
polynomial growth resolvent estimates are obtained, it is known from Lemma 4.7 in Burq \cite{Bu2004} that
these estimates, in fact, grow logarithmically (see also
 Tang and Zworski \cite{TZ00}). More precisely, our cited resolvent
$$
\chi \left(  \left. -\Delta \right\vert_D  - z^2 \right)^{-1} \chi 
$$ 
is a holomorphic $L^2(\Omega)$-valued function on $\{z\in\C:|\mathrm{Im}(z^2) | \leq \varepsilon_0\}$ and satisfies 
the assumptions in Lemma 4.7 from \cite{Bu2004}. Thus we conclude that there exist two constants $\lambda_0\gg 1$ and 
$0<\ep_0 \ll 1$ 
such that 
\begin{equation} \label{eq.main1}
    \lambda \big\| \chi \left(  \left. -\Delta \right\vert_D  - z^2 \right)^{-1} f \big\|_{L^2(\Omega)} 
    + \big\| \chi \nabla \left(  \left. -\Delta \right\vert_D  - z^2 \right)^{-1}  f \big\|_{L^2(\Omega)} 
    \leq C (\log \lambda)  \| f\|_{L^2(\Omega)}
\end{equation}
for $z^2=\lambda^2+\im \ep$ with $\lambda\geq \lambda_0$ and $0<\ep\leq \ep_0$, and 
for  any compactly supported functions $f \in L^2(\Omega)$. Hence we have the required estimate \eqref{eq.ME1}.
The proof of Theorem \ref{thm:main} is now finished.
\end{proof}

As a consequence of Theorem \ref{thm:main}, we have the following:
\begin{cor}\label{cor:main}
Let the assumptions in Theorem {\rm \ref{thm:main}} are fulfilled, and let $\lambda_0 \gg 1$ and $\ep_0 \ll 1$ be chosen 
as in \eqref{eq.main1} in the proof of Theorem {\rm \ref{thm:main}}.
 Then there exists a constant $C>0$ such that
\begin{equation*}
\sum_{k=0}^1 |\mu|^{1-k} \big\|\chi\nabla^k
\left(-\left. \Delta\right\vert_D - z^2\right)^{-1}  \chi
\big\|_{\mathcal{B}(L^2(\Omega))}\leq C \log(2+|\mu|)
\end{equation*}
for any $z= \mu+ \im \delta$
with $|\mu|>\mu_0$ and $\delta \in (0,\delta_0),$ where
\begin{equation} \label{eq.relation}
 \mu_0 \sim \lambda_0  \quad  \text{and} \quad \delta_0 \sim  \frac{\ep_0}{\lambda_0}. 
\end{equation}
\end{cor}
\begin{proof} Since
\[
z^2=\mu^2-\delta^2+2\im \mu \delta=\lambda^2 \pm \im \ep, 
\]
it follows that 
\[
|\mu|=\sqrt{\lambda^2+\delta^2} \quad \text{and} \quad 2|\mu| \delta=\ep,
\]
which implies the relation \eqref{eq.relation}.  
\end{proof}

We now turn to the proof of Theorem \ref{thm:mainB}.

\begin{proof}[Proof of Theorem {\rm \ref{thm:mainB}}]
Again we consider only the case when $\lambda>0$ is large.
The trivial estimate  $\log (2+\lambda) \lesssim (1+\lambda)^\alpha$ for any even integer  $\alpha>0$,  
the estimate of Theorem~\ref{thm:main} and Lemma~\ref{lem:lem 5.4}
imply that 
$$
\lambda \big\| \chi \left(  \left. -\Delta \right\vert_D  - z^2 \right)^{-1} f\big\|_{L^2(\Omega)} \lesssim
\| f\|_{H^\alpha(\Omega)}
$$ 
for any compactly supported $f \in H^\alpha(\Omega).$ Taking $\alpha=2$, we complete the proof of 
Theorem \ref{thm:mainB}.
\end{proof}

The following result is an immediate consequence of Theorem  \ref{thm:mainB}, and will be used in the proof of 
local energy decay. 

\begin{cor}\label{cor:main2}
Let the assumptions of Theorem {\rm \ref{thm:main}} are fulfilled, and let $\lambda_0 \gg 1$ and $\ep_0 \ll 1$ be chosen as in Corollary {\rm \ref{cor:main}}. 
 Then
\begin{equation*} 
|\mu| \big\| \chi \left(  \left. -\Delta \right\vert_D  - z^2 \right)^{-1} f\big\|_{L^2(\Omega)} \lesssim
\| f\|_{H^\alpha(\Omega)}
\end{equation*}
for some even integer $\alpha>0$, and for any $z= \mu+ \im \delta$
with $|\mu|>\mu_0$ and $\delta \in (0,\delta_0),$ where
\begin{equation*} 
\mu_0 \sim \lambda_0  \quad  \text{and} \quad \delta_0 \sim  \frac{\ep_0}{\lambda_0}, 
\end{equation*}
and for any  compactly supported $f \in H^\alpha(\Omega).$
\end{cor}


\section{Local energy decay: Proof of Theorem {\rm \ref{thmLED}}}
\label{sec:6}

We start with the following identity.
\begin{lem}
If $b>0$ and $t >0$, then
\begin{equation} \label{eq.ca1}
\int_{\R} \frac{e^{\im tz} z}{(1+z^2)(z^2-b^2)} \ dz = \frac{\im \cos (bt)}{2(1+b^2)}
- \frac{\pi \im e^{-t}}{ 1+b^2}.
\end{equation}
\end{lem}
\begin{proof}
Consider the contour
\[
 \Gamma_{R, \varepsilon} = \Gamma_R \cup (-R, -b-\varepsilon) \cup \gamma^-_\ep
\cup (-b+\ep, b-\ep) \cup \gamma^+_\ep \cup (b+\ep,R),
\]
where
\[
\Gamma_R = \{ z = R e^{\im \theta}: \theta \in [0,\pi] \},
\]
\[ \gamma^\pm_\ep = \{z = \pm b - \ep e^{-\im \theta}: \theta \in [0,\pi] \}.
\]
It is easy to check that
\begin{equation}\label{ca2}
    \left.\frac{e^{\im tz} z}{(1+z^2)(z^2-b^2)} \right|_{\Gamma_R} \longrightarrow 0
\end{equation}
as $ R \to \infty.$ Moreover, we can check that
\begin{equation}\label{ca3}
 \int_{\gamma^\pm_\ep} \frac{e^{\im tz} z}{(1+z^2)(z^2-b^2)} \, dz \longrightarrow
 -\frac{\pi \im}{2} \frac{e^{\pm \im tb}}{1+b^2}
\end{equation}
as $\ep \to 0.$
Applying the residue formula, we get
\[
\int_{\Gamma_{R,\ep}} \frac{e^{\im tz} z}{(1+z^2)(z^2-b^2)} \ dz
= -\pi \im \frac{e^{-t}}{1+b^2}.
\]
Thus, taking the limits of the above identity as $\ep\to0$ and $R\to\infty$, and
combining \eqref{ca2} and \eqref{ca3}, we get \eqref{eq.ca1}.
\end{proof}

We are now in a position to prove Theorem \ref{thmLED}.

\begin{proof}[Proof of Theorem {\rm \ref{thmLED}}] Take a cut-off function $\chi(x)$ such that $\chi(x)=1$ for $|x| \leq R$, such that the ball $|x| \leq R$ covers the obstacle.
For the sake of simplicity, we put $B = \sqrt{-\left. \Delta\right|_D}$.
The relation \eqref{eq.ca1} implies
\begin{equation} \label{eq.ca164}
\cos (tB)(1+B^2)^{-1} = -2\im \int_{\R} \frac{e^{\im tz} z}{1+z^2}
(z^2-B^2)^{-1} \ dz
+ 2 \pi  e^{-t} (1+B^2)^{-1}.
\end{equation}
Since the term
\[
 \chi \left(\int_{\R} \frac{e^{\im tz} z}{1+z^2}
(z^2-B^2)^{-1} \ dz \right) \chi
\]
is a $\mathcal{B}(H^2(\Omega),L^2(\Omega))$-valued meromorphic function having no  poles in the strip
\[
\{z\in\mathbb{C}: z = \mu  + \im \delta, \, \mu \in \mathbb{R}, \, 0 < \delta \ll 1\}
\]
 by Corollary \ref{cor:main2},
it is possible to shift the contour $\R$ into
a new one
\[
\gamma_c:=\{z=\mu+\im c:-\infty<\mu<\infty\} \quad \text{for $0<c \ll 1$},
\]
i.e., we have
\[
I_c(t):=\chi \left(\int_{\gamma_c} \frac{e^{\im tz} z}{1+z^2}
(z^2-B^2)^{-1} \ dz \right) \chi =\int_{\R} \frac{e^{\im tz} z}{1+z^2}
\chi (z^2-B^2)^{-1} \chi \, dz.
\]
By an elementary calculation $I_c(t)$ becomes
\[
I_c(t)=\int_\R \frac{e^{-ct} e^{\im t\mu}(\mu+\im c)}{1+(\mu+\im c)^2}
\chi ((\mu+\im c)^2-B^2)^{-1} \chi \, d\mu.
\]
Then applying \eqref{eq.ca164}, we find from the resolvent estimate \eqref{eq.ME1n1}
of Theorem \ref{thm:mainB}
 that
\begin{align*}
& \left\|\chi \cos \left(t B \right)(1+B^2)^{-1} \chi g \right\|_{L^2(\Omega)} \\
\lesssim &\,  e^{-ct} \int^\infty_0 \frac{|\mu+\im c|}{|1+(\mu+\im c)^2|} \left\| \chi ((\mu+\im c)^2-B^2)^{-1} \chi g
\right\|_{L^2(\Omega)} \, d\mu  \\
& \qquad \qquad +   e^{-t}\|\chi (1+B^2)^{-1} \chi g\|_{L^2(\Omega)}\\
\lesssim &\, e^{-ct} \left(\int^\infty_0 \frac{d\mu}{\sqrt{(\mu^2-c^2+1)^2+4c^2\mu^2}}
 \right) 
\|g \|_{H^2(\Omega)}
+   e^{-t}\|g\|_{L^2(\Omega)}\\
\lesssim &\, e^{-ct}\|g\|_{H^2(\Omega)}
\end{align*}
for any $t>0$. An interpolation between this estimate and the energy inequality gives \eqref{eq.led}.
The proof of Theorem \ref{thmLED} is complete.
\end{proof}


\section{Proof of Theorem \ref{thm:asymptotic expansion12}: Low frequency estimate}
\label{sec:7}

In this section we prove the resolvent estimate in low frequency which is stated in
Theorem \ref{thm:asymptotic expansion12}. Let $r_0>0$ be chosen so that
$$\Rn\setminus \Omega \subseteq \{x\in \Rn:|x|<r_0\}.
$$
Let us consider the free resolvent
\[
R_0(z):=\left(-\Delta-z^2\right)^{-1}
\quad \text{on $\Rn$}
\]
for $z\in \mathbb{C}$ with $|z|\leq \sigma_0$ for some $\sigma_0 \in (0,1]$.
The proof of Theorem \ref{thm:asymptotic expansion12} is based on
the following lemmas.
\begin{lem}\label{lem:Murata} Let $n\ge3$, $s>1/2$ and $\sigma_0 \in (0,1]$. Then
\begin{equation}\label{EQ:low expa}
\sum_{k=0}^1 |z|^{1-k} \left\|\langle x\rangle^{-s}\nabla^k
R_0(z)\langle x\rangle^{-s}
\right\|_{\mathcal{B}(L^2(\Rn))} \le C
\end{equation}
for any $z\in \mathbb{C}$ with $|z|\leq \sigma_0$.
\end{lem}

The estimate \eqref{EQ:low expa} is well known (see, e.g., Theorem 14.3.2 in \cite{Hormander2}), so we omit the proof.\\

If one needs to avoid the singularity of
$R_0(z)$ at $z=0$, we have the following:
\begin{lem}\label{lem:Murata1}
Let $n\ge3$, $s>1$ and $\sigma_0\in(0,1]$. Then
\begin{equation*}
\sum_{k=0}^1 \left\|\langle x\rangle^{-s}\nabla^k
R_0(z)\langle x\rangle^{-s}
\right\|_{\mathcal{B}(L^2(\Rn))} \leq C
\end{equation*}
for any $z\in \mathbb{C}$ with $|z|\leq \sigma_0$.
\end{lem}

The proof of Lemma \ref{lem:Murata1} follows from Lemma \ref{lem:A-3} in Appendix \ref{App:Appendix B}.  \\

We need the fact that zero is not eigenvalue of the Laplacian. Let us introduce weighted Sobolev spaces $H^m_s(\Omega)$
for a non-negative integer $m$ and real $s$:
\[
H^m_s(\Omega)=\{f: \langle x \rangle^s \pa^\alpha_x f\in L^2(\Omega) \}.
\]
In particular, we put
\[
L^2_s(\Omega)=H^0_s(\Omega).
\]
Then we have:
\begin{lem}\label{lem:Agmon}
Let $n\ge3$ and $s>1/2$.
Suppose that $u\in H^2_{\mathrm{loc}}(\Omega) \cap H^1_{-s}(\Omega)$ satisfies
\begin{equation*}
\left\{
\begin{aligned}
 \Delta u=0 \quad &\text{in $\Omega$,}\\
 u=0 \quad & \text{on $\pa\Omega$.}
 \end{aligned}
 \right.
\end{equation*}
Then $u(x)=0$ in $\Omega$.
\end{lem}
\begin{proof}
We claim that $u(x)$ is analytic in $\Omega$ and behaves like
\begin{equation}\label{EQ:ASYMPtotics}
\partial_x^\alpha u(x)=O(|x|^{-(n-2+|\alpha|)}), \quad |\alpha|\le1,
\end{equation}
as  $|x|\to\infty$.
To see this, we consider the extension of $u$ to $\Rn$.
Let $\psi(x)$ be a $C^\infty$ function on $\Rn$ such that
\[
\psi(x)=\left\{
\begin{aligned}
&0 & \quad (x\in U),\\
&1 &\quad (x\notin V),
\end{aligned}
\right.
\]
where  $U$ and $V$ are bounded and open neighbourhood of the obstacle $\mathcal{O}$ such that
\[
\mathcal{O} \subsetneq U \subsetneq V.
\]
We define $\widetilde{u}$ by letting
\[
\widetilde{u}(x)=\psi(x)u(x).
\]
We set
\begin{equation}\label{EQ:F}
f=-\Delta \widetilde{u}=-(\Delta \psi)u-2\nabla\psi \cdot \nabla u.
\end{equation}
Since $\widetilde{u} \in H^1_{-s}(\Rn)$ by assumption $u\in H^1_{-s}(\Omega)$, we can write
\begin{equation}\label{EQ:KERNEL}
\widetilde{u}(x)=\int_\Rn G_0(x-y) f(y)\, dy,
\end{equation}
where
\[
G_0(x)=\frac{\Gamma\left(\dfrac{n}{2}\right)}{4(n-2)\pi^{n/2}}\cdot \frac{1}{|x|^{n-2}}.
\]
On the other hand, since $u\in H^2_{\mathrm{loc}}(\Omega)$, it follows that $f\in H^1(\Rn)$.
Hence, by the regularity theory of elliptic equations, we see that $\widetilde{u}$ is $C^\infty$ on $\Rn$.
In particular, we deduce from \eqref{EQ:F} that $f$ is bounded and compactly supported on $\Rn$.
Since $\tilde{u}=u$ in $\Omega\setminus V$,
\eqref{EQ:ASYMPtotics} follow from \eqref{EQ:KERNEL}.

We prove that $u=0$ in $\Omega$. We put
\[
\Omega_R=\Omega\cap \{|x|<R\}
\]
for $R\gg 1$.
Integrating by parts, we have
\begin{equation} \label{EQ:PARTS}
\begin{split}
\int_{\Omega_R}|\nabla u(x)|^2\, dx
=&-\int_{\Omega_R} \Delta u(x) \overline{u(x)} \, dx
+\int_{\partial \Omega_R}
\frac{x}{|x|}\cdot \nabla u(x) \overline{u(x)} \, dS_R\\
=&\int_{|x|=R} \frac{x}{|x|}\cdot \nabla u(x) \overline{u(x)} \, dS_R.
\end{split}
\end{equation}
We observe from \eqref{EQ:ASYMPtotics} that
$$\frac{x}{|x|}\cdot \nabla u(x) \overline{u(x)}=O(|x|^{-2n+3})$$
as $|x|\to\infty$.
Then, letting $R\to \infty$ in \eqref{EQ:PARTS}, and noting $n\ge3$, we find that
\[
\int_{\Omega}|\nabla u(x)|^2\, dx=0,
\]
which implies that $u$ is constant in $\Omega$. Hence,
by using the boundary condition that $u=0$ on $\partial \Omega$,
we conclude that $u=0$ in $\Omega$.
The proof of Lemma \ref{lem:Agmon} is complete.
\end{proof}

We are now in a position to prove Theorem \ref{thm:asymptotic expansion12}.\\

\noindent
{\em Proof of Theorem {\rm \ref{thm:asymptotic expansion12}}.}
Let us introduce numbers $b$ and $d$ such that $d>b>r_0+3$ and fix them,
where
$r_0>0$ is chosen such that $\Rn\setminus \Omega \subseteq \{|x|\le r_0\}$.
Put $$\Omega_d=\Omega \cap \{|x|<d\}.$$
We consider the boundary value problem to the Poisson equation in the bounded
domain $\Omega_d$:
\begin{equation} \label{EQ:BDD}
\left\{
\begin{aligned}
-\Delta u=f \quad &\text{in $\Omega_d$},\\
u=0 \quad &\text{on $\partial\Omega_d$}.
\end{aligned}
\right.
\end{equation}
By the elliptic regularity theorem, for any $f\in L^2(\Omega_d)$,
there exists a unique solution $u\in H^2(\Omega_d)$ to \eqref{EQ:BDD}
such that
\[
\|u\|_{H^2(\Omega_d)}\le C\|f\|_{L^2(\Omega_d)}.
\]
Hence the mapping
of $f\in L^2(\Omega_d)$ to the unique solution
$u\in H^2(\Omega_d)$
determines an operator in $\mathcal{B}(L^2(\Omega_d),H^2(\Omega_d))$,
which is denoted by $L$.

Take $C^\infty$-functions $\varphi(x)$ and $\chi(x)$ such that
$\varphi(x)=1$ for $|x|\ge b$ and equal to $0$ for $|x|<b-1$;
$\chi(x)=1$ for $|x|\ge b-2$ and equal to $0$ for $|x|<b-3$.
For $f\in L^2_s(\Omega)$ with $s>1/2$, we define its restriction to $\Omega_d$
and zero extension as follows:
\[
f_d(x):=f(x)|_{\Omega_d} \quad \text{for $x\in \Omega_d$,}
\]
\[
f_0(x):=\left\{
\begin{aligned}
& f(x), &\quad x\in \Omega,\\
& 0,& \quad x\notin \Omega.
\end{aligned}
\right.
\]
Define an operator $R_1(z)$ by
\begin{equation}\label{EQ:RESOLVENT0}
R_1(z)f=\varphi R_0(z)(\chi f_0)+(1-\varphi)Lf_d, \quad f\in L^2_s(\Omega).
\end{equation}
Then we have $R_1(z)\in \mathcal{B}(L^2_s(\Omega),H^2_{-s}(\Omega))$ for $s>1/2$.
The operator thus defined obeys
\begin{equation}\label{EQ:RESOLVENT1}
(\left. -\Delta\right|_D-z^2)R_1(z)f=f+S(z)f \quad \text{in $\Omega$}
\end{equation}
for any $f\in L^2_s(\Omega)$, where $S(z)$ is defined by
\begin{equation}\label{EQ:RESOLVENT2}
S(z)f=-\left\{
2(\nabla \varphi)\cdot \nabla +\Delta \varphi
\right\}
\left\{ R_0(z)(\chi f_0)-Lf_d
\right\}-z^2(1-\varphi)Lf_d.
\end{equation}
It follows from Lemma \ref{lem:A-3} that $S(z)\in \mathcal{B}(L^2_s(\Omega))$
for $s>1/2$ and $z\in \mathbb{C}\setminus [0,\infty)$.
Since the support of $S(z)f$ is contained in $\overline{\Omega}_d$,
$S(z)$ is a compact operator in $L^2_s(\Omega)$.

\begin{lem} \label{lem:S-key-lemma}
Let $n\ge3$, and let $S(z)$ be the operator defined by \eqref{EQ:RESOLVENT2}. Then
the inverse $(I+S(z))^{-1}$ of $I+S(z)$ exists
as a $\mathcal{B}(L^2_s(\Omega))$-valued
meromorphic function of $z\in \mathbb{C}$
for $s>1/2$.
The set $\Lambda$ of poles is discrete and countable, and there exists an $\ep_0>0$ such that
$\Lambda$ has no intersection
with $\{z\in \mathbb{C}:|z|<\varepsilon_0\}$.
In addition, $(I+S(z))^{-1}$ has the same type of estimate as
\eqref{EQ:low expa} from Lemma {\rm \ref{lem:Murata}}.
More precisely, the operator $(I+S(z))^{-1}$
is obtained as a Neumann series expansion{\rm :} For
$|z|<\varepsilon_0$,
\begin{equation}\label{EQ:Neumann}
(I+S(z))^{-1}=(I+S(0))^{-1}\sum_{j=0}^\infty
\left[(S(0)-S(z))(I+S(0))^{-1} \right]^j.
\end{equation}
\end{lem}
\begin{proof}
To begin with, we claim that
\[
(I+S(0))^{-1}\in \mathcal{B}(L^2_s(\Omega)).
\]
Indeed, if we prove that $I+S(0)$ is injective, the conclusion follows
from Fredholm's alternative, since the operator $S(0)$ is compact.
Therefore, for the time being,
we concentrate on proving the injectivity of $I+S(0)$.
Let
us assume that
\[
(I+S(0))f=0, \quad f\in L^2_s(\Omega).
\]
Then it follows from \eqref{EQ:RESOLVENT1} that $R_1(0)f$ satisfies the elliptic boundary value problem
of the form
\begin{equation*}
\left\{
\begin{aligned}
 \Delta R_1(0)f=0 & \quad \text{in $\Omega$}.\\
 R_1(0)f=0 &\quad \text{on $\pa\Omega$.}
\end{aligned}
\right.
\end{equation*}
Since
\[
R_1(0)f \in H^2_{\mathrm{loc}}(\Omega) \cap H^1_{-s}(\Omega),
\]
we deduce from Lemma \ref{lem:Agmon} that
\begin{equation}\label{EQ:ZERO}
R_1(0)f=0 \quad \text{in $\Omega$.}
\end{equation}
Since
\[
\text{$R_1(0)f= R_0(0)(\chi f_0)$ \quad for $|x|\ge b$}
\]
by \eqref{EQ:RESOLVENT0}, the identity \eqref{EQ:ZERO} together with \eqref{EQ:RESOLVENT0} imply that
\begin{equation}\label{EQ:TRIVIALITY}
R_0(0)(\chi f_0)=0 \quad \text{for $|x|\ge b$,}
\end{equation}
\begin{equation}\label{EQ:TRIVIALITY1}
Lf_d=0 \quad \text{for $x \in \Omega$, $|x|\le b-1$.}
\end{equation}
Since we have the relation $-\Delta R_0(0)(\chi f_0)=\chi f_0$ in $\Rn$, it follows from
\eqref{EQ:TRIVIALITY} that $\chi f_0=0$ for $|x|\ge b$, i.e.,
\begin{equation}\label{EQ:0-ZERO}
f(x)=0 \quad \text{for $|x|\ge b$.}
\end{equation}
Similarly, $Lf_d$ satisfies the equation $-\Delta Lf_d=f_d$ in
$\Omega_d$, and hence, by using \eqref{EQ:TRIVIALITY1}, we get
\[
f(x)=0\quad \text{for $x \in \Omega$, $|x|\le b-1$.}
\]
These imply that $\chi f_0=f_0$ and
\begin{equation} \label{EQ:Poisson1}
\left\{
\begin{aligned}
-\Delta R_0(0)(\chi f_0)&=f_0 &\quad \text{in $\Rn$,}\\
R_0(0)(\chi f_0)&=0 &\quad \text{on $|x|=d$.}
\end{aligned}
\right.
\end{equation}
On the other hand, if we define
\[
v=\begin{cases}
Lf_d \quad &\text{in $\Omega_d$,}\\
0 \quad & \text{in $\Rn\setminus \Omega$,}
\end{cases}
\]
then we see from the elliptic regularity theorem that
$v\in H^2(B_d(0))$, and
\begin{equation} \label{EQ:Poisson2}
\left\{
\begin{aligned}
-\Delta v&=f_0 \quad &\text{in $B_d(0)$,}\\
v&=0 \quad &\text{on $|x|=d$.}
\end{aligned}
\right.
\end{equation}
Hence it follows from \eqref{EQ:Poisson1} and \eqref{EQ:Poisson2} that
$R_0(0)(\chi f_0)=v$ in $B_d(0)$, and hence,
\[
R_0(0)(\chi f_0)=Lf_d \quad \text{in $\Omega_d$,}
\]
which implies that
$R_1(0)f=Lf_d$ in $\Omega_d$.
By this relation and \eqref{EQ:ZERO} we have
\[
0=-\Delta R_1(0)f=-\Delta Lf_d=f_d \quad \text{in $\Omega_d$,}
\]
i.e., $f=0$ in $\Omega_d$, which together with \eqref{EQ:0-ZERO} shows
$f=0$ in $\Omega$. This proves the injectivity of $I+S(0)$.

Put $M=\|(I+S(0))^{-1}\|_{\mathcal{B}(L^2_s(\Omega))}$.
By the continuity of $S(z)$, there exists $\varepsilon_0>0$
such that
\[
\|S(0)-S(z)\|_{\mathcal{B}(L^2_s(\Omega))}<\frac{1}{2M}
\]
for any $z\in \{ z\in \mathbb{C}:|z|<\varepsilon_0\}$. Thus the inverse $(I+S(z))^{-1}$
is obtained as a Neumann series expansion: For
$|z|<\varepsilon_0$,
\begin{equation*}
(I+S(z))^{-1}=(I+S(0))^{-1}\sum_{j=0}^\infty
\left[(S(0)-S(z))(I+S(0))^{-1} \right]^j,
\end{equation*}
which proves \eqref{EQ:Neumann}.
Since $S(z)$ is holomorphic in $\{z\in \mathbb{C}:|z|<\ep_0\}$, applying analytic
Fredholm's alternative, we conclude from
\cite[p.~592, Lemma~13]{Dunford} that $(I+S(z))^{-1}$ exists
in $\{z\in \mathbb{C}:|z|<\ep_0\}$
as a meromorphic function, and the set $\Lambda$
of the poles is discrete and countable in $\{z\in \mathbb{C}:|z|<\ep_0\}$.
Thus Lemma \ref{lem:S-key-lemma} follows from
Lemma \ref{lem:Murata} and \eqref{EQ:Neumann}. The proof of Lemma \ref{lem:S-key-lemma}
is complete.
\end{proof}

We are now in a position to complete the proof of Theorem \ref{thm:asymptotic expansion12}.
\begin{proof}[Completion of the proof of
Theorem {\rm \ref{thm:asymptotic expansion12}}]
Define
\begin{equation} \label{EQ:res}
\widetilde{R}(z)=R_1(z)(I+S(z))^{-1}.
\end{equation}
Then we conclude
from Lemma \ref{lem:Murata}, Lemma
\ref{lem:S-key-lemma} and \eqref{EQ:res} that
there exists $\varepsilon_0>0$
such that $\{z\in \mathbb{C}:|z|<\ep_0\}\cap \Lambda=\emptyset$ and
one has the estimate
\begin{equation*}
\sum_{k=0}^1 |z|^{1-k}\big\| \langle x \rangle^{-s}\nabla^k
\widetilde{R}(z)\langle x \rangle^{-s}\big\|_{\mathcal{B}(L^2(\Omega))}\le C
\end{equation*}
for $z\in \{z\in \mathbb{C}:|z|<\ep_0\}$ and $s>1/2$. Furthermore,
since the resolvent set of $-\Delta|_D$ contains $\{z\in \mathbb{C}:|z|<\ep_0\}$,
it follows from \eqref{EQ:RESOLVENT1} and \eqref{EQ:res} that the identity
\[
\langle x \rangle^{-s}\widetilde{R}(z)\langle x \rangle^{-s}
=\langle x \rangle^{-s}(-\Delta|_D-z^2)^{-1}\langle x \rangle^{-s} \quad \text{on $L^2(\Omega)$}
\]
holds for all $z\in \{z\in \mathbb{C}:|z|<\ep_0\}$.
The proof of Theorem \ref{thm:asymptotic expansion12} is now finished.
\end{proof}


\appendix

\section{Some examples of domains}
\label{Append:Appendix A}

In this appendix we shall provide some examples of obstacles satisfying the assumptions
or not satisfying.
\begin{exam}\label{exam1}
{\rm First we consider the case of two convex obstacles in $\mathbb{R}^3$ which is the classical example due to 
Ikawa \cite{Ikawa}.
In this case the $x$-projection of the trapped set is a segment, while a $\delta$-neighbourhood $\Sigma_\delta$ of the trapped set can be covered by
\[
  \bigcup_{j=1}^{N(\delta)} (K_{j,\delta} \times E_{j,\delta}), 
\]
where $N(\delta) \sim 1/\delta$, so Assumption B is true. 
The sojourn time is  
\[
t_{\delta,R,\rho} \leq C(R) \log \left( \frac{1}{\delta}\right),
\]
so Assumption C is fulfilled too.
}
\end{exam}

\begin{figure}[H]
\includegraphics[width=0.85\textwidth]{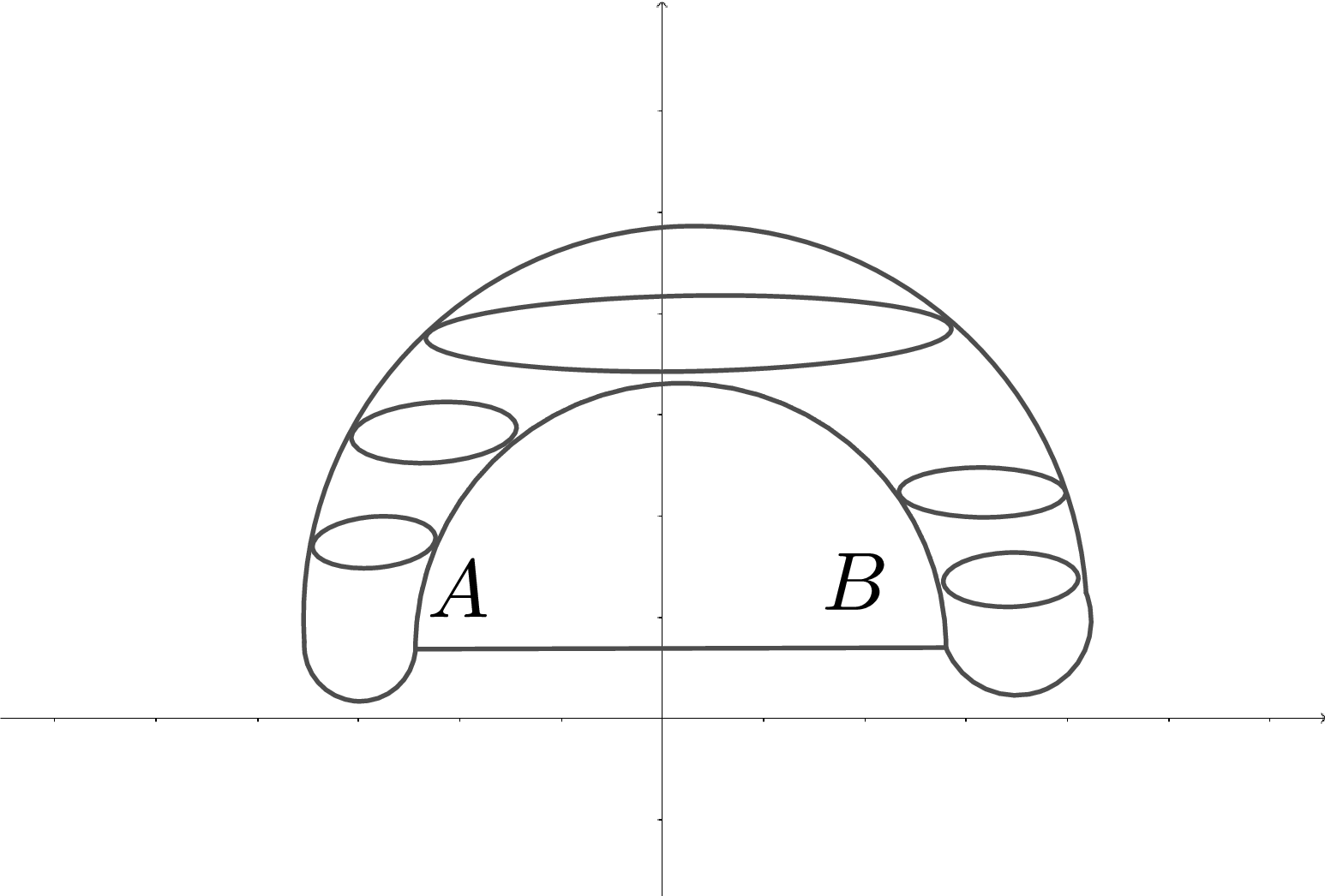}
\caption{Trapped exterior domain satisfying Assumptions B and C.}
\label{fig:Pic1}
\end{figure}
\begin{exam}\label{exam2}
{\rm Another case is shown on Figure \ref{fig:Pic1}. It is clear that we have the situation similar to that of Example \ref{exam1} and we have a unique trapped segment $AB.$ In this case Assumptions B and C are satisfied.
$\Sigma_\delta$ is covered by
\[
  \bigcup_{j=1}^{N(\delta)} (K_{j,\delta} \times E_{j,\delta}), 
\]
where $N(\delta) \sim 1/\delta$, so Assumption B is true. 
The sojourn time is  
\begin{equation}\label{stje1}
    t_{\delta,R,\rho} \leq C(R) \log \left( \frac{1}{\delta}\right),
\end{equation}
so Assumption C is fulfilled too. 
}
\end{exam}

\begin{exam}\label{exam3} 
{\rm If we rotate the obstacle on Figure \ref{fig:Pic1} around the vertical axis, then we can consider an example in $\mathbb{R}^3$, and we shall have a disk in the place of the segment as $x$-projection of the trapped set. It is easy to see that
$ \Pi^{-1}(x)$
can be covered by $N(\delta)\sim \delta^{-2}$ disks of radius $\delta$ on $\mathbb{S}^2$, so Assumption B is not fulfilled. One can try to verify \eqref{stje1}, and we can immediately say that our main result is not applicable in this case, since Assumption B is not fulfilled.
}
\end{exam}
\begin{figure}[H]
\includegraphics[width=0.85\textwidth]{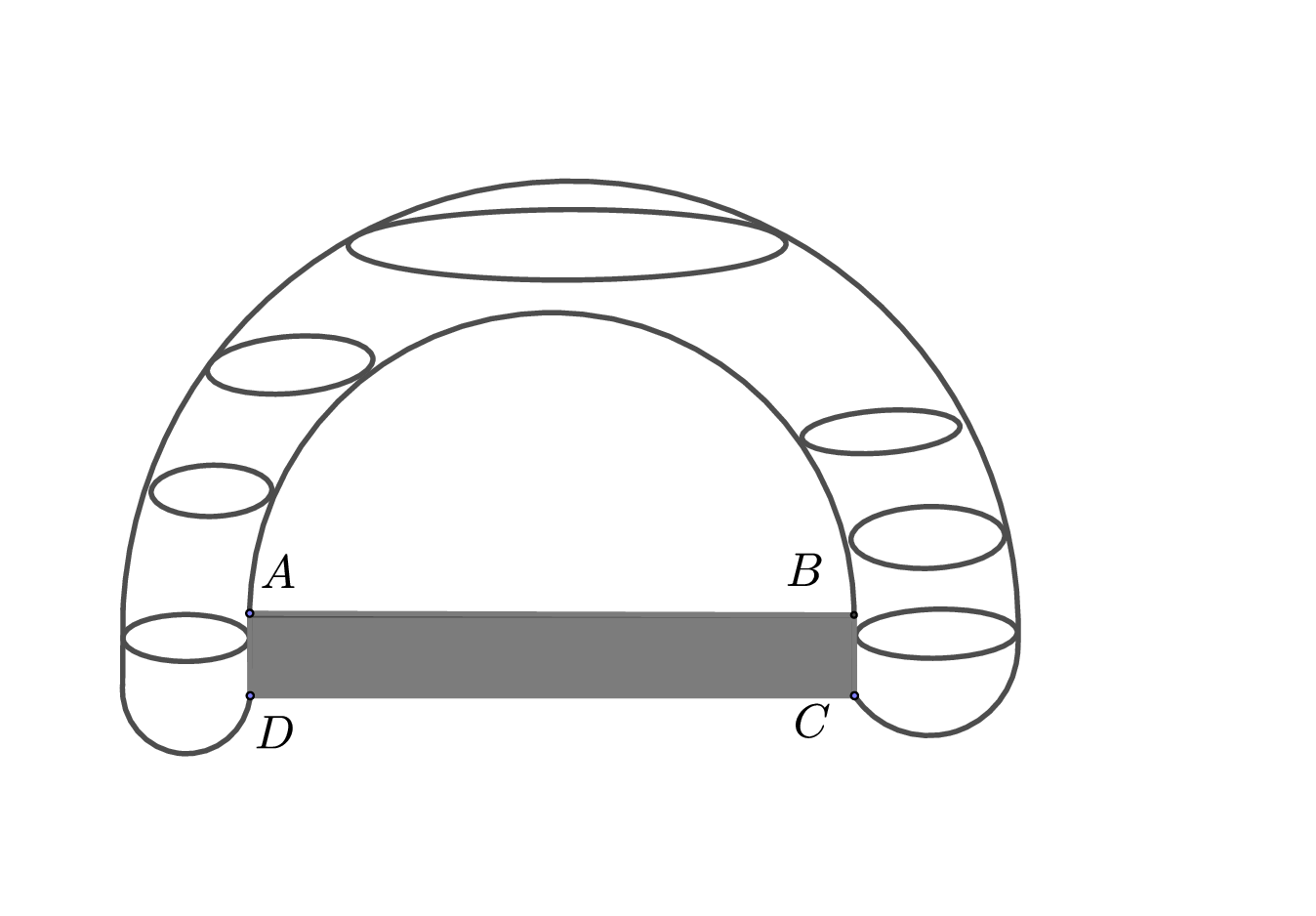}
\caption{Trapped exterior domain satisfying Assumptions B, but not C.}
\label{fig:Pic2}
\end{figure}

\begin{exam}\label{exam4} 
{\rm Figure \ref{fig:Pic2} represents an obstacle with $x$-projection of the trapped set given by the rectangle $ABCD.$ This example is a modification of Example \ref{exam2} and one can check that Assumption B is satisfied, while Assumption C is not satisfied, since
\begin{equation*} 
    t_{\delta,R,\rho} \leq \frac{C(R)}{\delta}.
\end{equation*}
}
\end{exam}


\section{(Commutator estimates)}
\label{App:Appendix A}
In this appendix we prove the estimate \eqref{eq.comm} in Section \ref{sec:4}. The following formula is well known.
\begin{lem} \label{lem:well}
Let $A$ be a non-negative self-adjoint operator on a Hilbert space.
Then
\begin{equation*} 
(A - z)^{- \beta} =
\frac{1}{\Gamma(\beta)}
\int^{\infty}_{0} t^{\beta - 1} e^{z t}
	e^{-t A} \,dt
\end{equation*}
for any $z \in \mathbb{C}$ with
${\rm Re}(z) < \inf{\sigma(A)}$ and
$\beta > 0$, where $\Gamma(\beta)$ stands for the Gamma function, and
$\sigma(A)$ denotes the spectrum of $A$.
\end{lem}

We prepare the commutator estimate for the heat semi-group $\{e^{\left. t \Delta \right|_D}\}$.

\begin{lem} \label{lem:lem A}
Assume that $0<s\leq1$. Let $f\in H^{N-1}(\Omega)$ for some non-negative integer $N$. Then
\begin{equation*} 
\left\|\left[\langle x \rangle^s ,e^{\left. t \Delta \right|_D }\right]f\right\|_{H^N(\Omega)}
\lesssim \|f \|_{H^{N-1}(\Omega)}
\end{equation*}
for any $t>0$. 
\end{lem}
\begin{proof} 
Let $u(t,x):=(e^{t \left. \Delta\right|_D} f)(x)$ be a solution to the initial-boundary value problem
for the heat equation with data $f$:
\begin{equation*} 
\left\{
\begin{array}{ll}
(\pa_t-\left. \Delta\right|_D )u=0, & t>0, \quad x\in \Omega,\\
 u(0,x)=f, &  x\in\Omega. 
 \end{array} \right.
\end{equation*}
Then the function
\[
u_s(t,x):=\langle x \rangle^s u(t,x)
\]
solves the following initial-boundary value problem:
\begin{equation*} 
\left\{
\begin{array}{ll}
(\pa_t-\left. \Delta\right|_D)u_s=\left[\left. \Delta\right|_D,\langle x \rangle^s \right]u, & t>0, \quad x\in \Omega,\\
 u_s(0,x)=\langle x \rangle^s f, &  x\in\Omega. 
\end{array} \right.
\end{equation*}
If $v(t,x):=e^{t\left. \Delta\right|_D} \langle x \rangle^s f$, then $v$ solves the following homogeneous problem:
\begin{equation*} 
\left\{
\begin{array}{ll}
(\pa_t-\left. \Delta\right|_D)v=0, & t>0, \quad x\in \Omega,\\
 v(0,x)=\langle x \rangle^s f, &  x\in\Omega. 
\end{array} \right.
\end{equation*}
We note that  $u_s(t)-v(t)$ satisfies
\[
u_s(t)-v(t)=\left[\langle x \rangle^s,e^{t \left. \Delta\right|_D} \right] f,
\]
and solves the problem
\begin{equation*} 
\left\{
\begin{array}{ll}
(\pa_t-\left. \Delta\right|_D)(u_s-v)=\left[\left. \Delta\right|_D,\langle x \rangle^s \right]u, & t>0, \quad x\in \Omega,\\
(u_s-v)(0,x)=0, &  x\in\Omega. 
\end{array} \right.
\end{equation*}
Then applying Theorem 5 in p. 360 from Evans \cite{Evans},
we estimate
\[
\|u_s-v\|_{L^2_t H^{N+1}(\Omega)}+\|u_s-v\|_{L^\infty_t H^N(\Omega)}
\lesssim \left\| \left[ \left. \Delta\right|_D,\langle x \rangle^s \right] u \right\|_{L^2_tH^{N-1}(\Omega)}
\lesssim \|u\|_{L^2_t H^N(\Omega)}
\]
for any $N\geq 0$, provided that $0<s\leq 1$.
Finally, we use a standard smoothing estimate for the heat equation:
\[
 \|u\|_{L^2_t H^{N}(\Omega)}\leq  \|f\|_{H^{N-1}(\Omega)} \quad
 \text{for any $N \geq 0$.}
\]
This estimate is very well known for $N \geq 1,$ but using the fact that $(1-\left. \Delta\right|_D)^{-1/2}$ commutes with the operator $\pa_t-\left. \Delta \right|_D$, we see that this estimate is  true in the case $N=0.$
Hence, we can conclude that
\[
\left\|\left[\langle x \rangle^s,e^{t \left. \Delta\right|_D} \right]f \right\|_{H^N}
=\| u_s-v \|_{H^N}\lesssim \|f\|_{H^{N-1}(\Omega)}
\]
for any $t>0$ and $N\geq 0$. The proof of Lemma \ref{lem:lem A} is complete.
\end{proof}

Based on Lemmas \ref{lem:well} and \ref{lem:lem A}, we prove the commutator estimate.
\begin{lem} \label{lem:lem A2}
Let $0<s\leq1$. and let $N$ be an integer with $N\geq 0$. Then
\begin{equation}\label{eq.k-comm}
\left\| \left[\langle x \rangle^s ,(1-\left. \Delta \right|_D)^{-\frac{k}{2}}\right]f\right\|_{H^N(\Omega)}
\lesssim \| f \|_{H^{N-k}(\Omega)}
\end{equation}
for any $f \in H^N(\Omega)$ and integer $k$ with $1\leq k\leq N$. When $N=0$,
\eqref{eq.k-comm} is valid for $k=1$.
\end{lem}
\begin{proof}
Applying Lemma \ref{lem:well}, we have the following identity:
\begin{equation}\label{eq.ind}
\left[\langle x \rangle^s ,(1-\left. \Delta \right|_D)^{-\frac{k}{2}}\right]f
=\frac{1}{\Gamma\left(\frac{k}{2}\right)}
\int^\infty_0 t^{\frac{k}{2}-1} e^{-t} \left[\langle x \rangle^s ,e^{t\left. \Delta \right|_D} \right]f\, dt
\end{equation}
for any $k>0$.
Hence, it follows from Lemma \ref{lem:lem A} that if $k=1$, then
\begin{equation}\label{eq:1-ind}
\begin{split}
\left\| \left[\langle x \rangle^s ,(1-\left. \Delta \right|_D)^{-\frac{1}{2}}\right]f \right\|_{H^N(\Omega)}
\leq & \, \frac{1}{\Gamma\left(\frac{1}{2}\right)}
\int^\infty_0 t^{-\frac{1}{2}} e^{-t} \left\| \left[\langle x \rangle^s ,e^{t\left. \Delta \right|_D} \right]f \right\|_{H^N(\Omega)} \, dt\\
\lesssim &\, \|f\|_{H^{N-1}(\Omega)}
\end{split}
\end{equation}
for $N\geq 0$. Therefore, \eqref{eq.k-comm} is true for $(N,k)=(0,1)$, $(1,1)$.

We must prove the case when $N>1$.
Based on the identity \eqref{eq.ind},
we prove \eqref{eq.k-comm} by induction argument on $k$. Thanks to \eqref{eq:1-ind},
\eqref{eq.k-comm} is true for $k=1$.
Suppose that \eqref{eq.k-comm} is true for $k=N-1$, i.e.,
\[
\left\| \left[\langle x \rangle^s ,(1-\left. \Delta \right|_D)^{-\frac{N-1}{2}}\right] f \right\|_{H^N(\Omega)}
\lesssim \|f\|_{H^1(\Omega)},
\]
which is equivalent to the following:
\begin{equation}\label{eq.N-ind}
\left\| (1-\left. \Delta \right|_D)^{\frac{N}{2}}
\left[\langle x \rangle^s ,(1-\left. \Delta \right|_D)^{-\frac{N-1}{2}}\right] (1-\left. \Delta \right|_D)^{-\frac{1}{2}}f \right\|_{L^2(\Omega)}
\lesssim \|f\|_{L^2(\Omega)}.
\end{equation}
We note that
\begin{equation}\label{eq.Int-comm}
\begin{split}
& (1-\left. \Delta \right|_D)^{\frac{N}{2}}\left[\langle x \rangle^s ,(1-\left. \Delta \right|_D)^{-\frac{N}{2}}\right]\\
=&\, (1-\left. \Delta \right|_D)^{\frac{N}{2}}
\left[\langle x \rangle^s ,(1-\left. \Delta \right|_D)^{-\frac{N-1}{2}}\right] (1-\left. \Delta \right|_D)^{-\frac{1}{2}}\\
& +(1-\left. \Delta \right|_D)^{\frac{1}{2}}\left[\langle x \rangle^s,(1-\left. \Delta \right|_D)^{-\frac12} \right]
\end{split}
\end{equation}
Indeed, since
\begin{align*}
\left[\langle x \rangle^s ,(1-\left. \Delta \right|_D)^{-\frac{N}{2}}\right]
=\langle x \rangle^s (1-\left. \Delta \right|_D)^{-\frac{N-1}{2}}(1-\left. \Delta \right|_D)^{-\frac{1}{2}}
-(1-\left. \Delta \right|_D)^{-\frac{N}{2}}\langle x \rangle^s,
\end{align*}
it follows that
\begin{equation} \label{eq.arrow1}
\begin{split}
& (1-\left. \Delta \right|_D)^{\frac{N}{2}}\left[\langle x \rangle^s ,(1-\left. \Delta \right|_D)^{-\frac{N}{2}}\right]\\
=&\, (1-\left. \Delta \right|_D)^{\frac{N}{2}}\langle x \rangle^s (1-\left. \Delta \right|_D)^{-\frac{N-1}{2}}(1-\left. \Delta \right|_D)^{-\frac{1}{2}}
-\langle x \rangle^s.
\end{split}
\end{equation}
Here, the first term in the right becomes
\begin{equation}\label{eq.arrow2}
\begin{split}
&(1-\left. \Delta \right|_D)^{\frac{N}{2}}\langle x \rangle^s (1-\left. \Delta \right|_D)^{-\frac{N-1}{2}}(1-\left. \Delta \right|_D)^{-\frac{1}{2}}\\
=&\, (1-\left. \Delta \right|_D)^{\frac{N}{2}} \left[\langle x \rangle^s ,(1-\left. \Delta \right|_D)^{-\frac{N-1}{2}}\right]
(1-\left. \Delta \right|_D)^{-\frac{1}{2}} \\
& \, +(1-\left. \Delta \right|_D)^{\frac{1}{2}}\langle x \rangle^s(1-\left. \Delta \right|_D)^{-\frac{1}{2}}\\
=&\, (1-\left. \Delta \right|_D)^{\frac{N}{2}} \left[\langle x \rangle^s ,(1-\left. \Delta \right|_D)^{-\frac{N-1}{2}}\right]
(1-\left. \Delta \right|_D)^{-\frac{1}{2}}\\
&\, +(1-\left. \Delta \right|_D)^{\frac{1}{2}}\left[\langle x \rangle^s ,(1-\left. \Delta \right|_D)^{-\frac{1}{2}}\right]+\langle x \rangle^s .
\end{split}
\end{equation}
Hence, \eqref{eq.Int-comm} follows from \eqref{eq.arrow1} and \eqref{eq.arrow2}.
Thus, combining \eqref{eq.N-ind} and \eqref{eq.Int-comm}, we conclude that
\[
\left\|\left[\langle x \rangle^s ,(1-\left. \Delta \right|_D)^{-\frac{N}{2}}\right]f \right\|_{H^N(\Omega)}
\lesssim \|f\|_{L^2(\Omega)}
+\left\|\left[\langle x \rangle^s,(1-\left. \Delta \right|_D)^{-\frac12} \right] f \right\|_{H^1(\Omega)}
\lesssim \|f\|_{L^2(\Omega)},
\]
where we used \eqref{eq:1-ind} for $N=1$ in the last step. This proves \eqref{eq.k-comm} for $k=N$.
In conclusion, \eqref{eq.k-comm} is true for any integer $k$ with $1\leq k\leq N$. The proof of Lemma \ref{lem:lem A2} is complete.
\end{proof}

We are now in a position to prove estimate \eqref{eq.comm}. 
\begin{lem} \label{lem:A5}
Let $\Omega^\prime$ be an unbounded subdomain of $\Omega$ with a smooth boundary.
Let $N$ be a positive integer, $w \in H^N(\Omega)$ and $0<s\leq1$. Then one has the following estimate{\rm :}
\begin{equation} \label{eq.commutator-0}
 \sum_{|\alpha| \leq N}\|\langle x \rangle^{-s} \partial_x^\alpha w \|_{L^2(\Omega^\prime)} \sim \Big\|\langle x \rangle^{-s} \left( 1  \left. -\Delta \right\vert_D \right)^{\frac{N}{2}} w \Big\|_{L^2(\Omega^\prime)}.
 \end{equation}
\end{lem}
\begin{proof}
If we prove that
\begin{equation} \label{eq.k-star}
\Big\| \left( 1 \left. -\Delta \right\vert_D \right)^{\frac{k}{2}} \langle x \rangle^{-s} w\Big\|_{L^2(\Omega^\prime)}
\sim
 \Big\|\langle x \rangle^{-s} \left( 1  \left. -\Delta \right\vert_D \right)^{\frac{k}{2}} w \Big\|_{L^2(\Omega^\prime)}
\end{equation}
for $0\leq k\leq N$, then, by combining the elliptic regularity theorem, we get the required estimate \eqref{eq.commutator-0}.
For the proof of \eqref{eq.k-star}, we have only to show that
\begin{equation} \label{eq.star}
\Big\| \left( 1 \left. -\Delta \right\vert_D \right)^{\frac{1}{2}} \langle x \rangle^{-s} w\Big\|_{L^2(\Omega^\prime)}
\sim
 \Big\|\langle x \rangle^{-s} \left( 1  \left. -\Delta \right\vert_D \right)^{\frac{1}{2}} w \Big\|_{L^2(\Omega^\prime)},
\end{equation}
since the proof of \eqref{eq.k-star} for $2\leq k\leq N$ can be performed in completely similar way as that of \eqref{eq.star}.

To begin with, we prove that
\begin{equation}\label{eq.with}
 \Big\|\langle x \rangle^{-s} \left( 1  \left. -\Delta \right\vert_D \right)^{\frac{1}{2}} w \Big\|_{L^2(\Omega^\prime)}
\lesssim \Big\| \left( 1 \left. -\Delta \right\vert_D \right)^{\frac{1}{2}} \langle x \rangle^{-s} w\Big\|_{L^2(\Omega^\prime)}.
\end{equation}
If we put
\[
f:=\left( 1 \left. -\Delta \right\vert_D \right)^{\frac{1}{2}} \langle x \rangle^{-s} w,
\]
then \eqref{eq.with} is equivalent to the following:
\[
 \Big\|\langle x \rangle^{-s} \left( 1  \left. -\Delta \right\vert_D \right)^{\frac{1}{2}} \langle x \rangle^s
 \left( 1  \left. -\Delta \right\vert_D \right)^{-\frac{1}{2}}f \Big\|_{L^2(\Omega^\prime)}
\lesssim \|f\|_{L^2(\Omega^\prime)}.
\]
Hence, by duality, we have only to show that
\begin{equation}\label{eq.dual}
 \Big\| \left( 1  \left. -\Delta \right\vert_D \right)^{-\frac{1}{2}}\langle x \rangle^s \left( 1  \left. -\Delta \right\vert_D \right)^{\frac{1}{2}}
\langle x \rangle^{-s}g \Big\|_{L^2(\Omega^\prime)}
\lesssim \|g\|_{L^2(\Omega^\prime)}.
\end{equation}
Writing
\[
\left( 1  \left. -\Delta \right\vert_D \right)^{-\frac{1}{2}}\langle x \rangle^s \left( 1  \left. -\Delta \right\vert_D \right)^{\frac{1}{2}}
\langle x \rangle^{-s} =I+[\left( 1 \left. -\Delta \right\vert_D \right)^{-\frac{1}{2}},\langle x \rangle^s]
\left( 1 \left. -\Delta \right\vert_D \right)^{\frac{1}{2}}\langle x \rangle^{-s},
\]
we use Lemma \ref{lem:lem A2} to deduce that
\begin{align*}
 & \Big\| \left( 1  \left. -\Delta \right\vert_D \right)^{-\frac{1}{2}}\langle x \rangle^s \left( 1  \left. -\Delta \right\vert_D \right)^{\frac{1}{2}}
\langle x \rangle^{-s}g \Big\|_{L^2(\Omega^\prime)}\\
\lesssim &\, \|g\|_{L^2(\Omega^\prime)}
+\Big\| \left( 1  \left. -\Delta \right\vert_D \right)^{-\frac{1}{2}}\left( 1  \left. -\Delta \right\vert_D \right)^{\frac{1}{2}}
\langle x \rangle^{-s} g \Big\|_{L^2(\Omega^\prime)} \\
\leq &\, 2\|g\|_{L^2(\Omega^\prime)},
\end{align*}
which proves \eqref{eq.dual}.

Next, we prove the converse, i.e.,
\begin{equation}\label{eq:converse}
\Big\| \left( 1 \left. -\Delta \right\vert_D \right)^{\frac{1}{2}} \langle x \rangle^{-s} w\Big\|_{L^2(\Omega^\prime)}
\lesssim
 \Big\|\langle x \rangle^{-s} \left( 1  \left. -\Delta \right\vert_D \right)^{\frac{1}{2}} w \Big\|_{L^2(\Omega^\prime)}.
\end{equation}
If we put
\[
h:=\langle x \rangle^{-s} \left( 1  \left. -\Delta \right\vert_D \right)^{\frac{1}{2}} w,
\]
then \eqref{eq:converse} is equivalent to the following:
\begin{equation*} 
\Big\|
\left( 1  \left. -\Delta \right\vert_D \right)^{\frac{1}{2}}\langle x \rangle^{-s}
\left( 1  \left. -\Delta \right\vert_D \right)^{-\frac{1}{2}}\langle x \rangle^{s}h
\Big\|_{L^2(\Omega^\prime)}\lesssim \|h\|_{L^2(\Omega^\prime)}.
\end{equation*}
By duality, we have only to show that
\begin{equation} \label{eq.last}
\Big\|
\langle x \rangle^s \left( 1  \left. -\Delta \right\vert_D \right)^{-\frac{1}{2}}
\langle x \rangle^{-s}\left( 1  \left. -\Delta \right\vert_D \right)^{\frac{1}{2}} \va
\Big\|_{L^2(\Omega^\prime)}\lesssim \|\va\|_{L^2(\Omega^\prime)}.
\end{equation}
Since
\[
\langle x \rangle^s \left( 1  \left. -\Delta \right\vert_D \right)^{-\frac{1}{2}}
\langle x \rangle^{-s}\left( 1  \left. -\Delta \right\vert_D \right)^{\frac{1}{2}}
=I+\left[\langle x \rangle^s ,\left( 1  \left. -\Delta \right\vert_D \right)^{-\frac{1}{2}} \right]
\langle x \rangle^{-s}\left( 1  \left. -\Delta \right\vert_D \right)^{\frac{1}{2}},
\]
we estimate, by using Lemma \ref{lem:lem A2} twice,
\begin{align*}
& \Big\|
\langle x \rangle^s \left( 1  \left. -\Delta \right\vert_D \right)^{-\frac{1}{2}}
\langle x \rangle^{-s}\left( 1  \left. -\Delta \right\vert_D \right)^{\frac{1}{2}}\va
\Big\|_{L^2(\Omega^\prime)}\\
\lesssim & \, \|\va \|_{L^2(\Omega^\prime)}+\Big\|\left( 1  \left. -\Delta \right\vert_D \right)^{-\frac{1}{2}} \langle x \rangle^{-s}\left( 1  \left. -\Delta \right\vert_D \right)^{\frac{1}{2}}\va \Big\|_{L^2(\Omega^\prime)}\\
\lesssim & \, \|\va \|_{L^2(\Omega^\prime)}+\Big\|\left( 1  \left. -\Delta \right\vert_D \right)^{-\frac{1}{2}}
\left( 1  \left. -\Delta \right\vert_D \right)^{\frac{1}{2}}\va \Big\|_{L^2(\Omega^\prime)}\\
=&\, 2\|\va \|_{L^2(\Omega^\prime)}.
\end{align*}
This proves \eqref{eq.last}.

In conclusion, \eqref{eq.star} follows from \eqref{eq.with} and \eqref{eq:converse}.
The proof of Lemma \ref{lem:A5} is finished.
\end{proof}


 \section{(Free resolvent estimates)}
 \label{App:Appendix B}

The next lemma in the particular case $n=3$ can be found in \cite{Georg-Visciglia}
and for general case $n \geq 3$ in \cite{Comech}, where an idea from \cite{Ginibre} is used.
\begin{lem} \label{lem:A-3}
Let $n\ge3$, $s_1 >1/2$, $s_2 >1/2$ and $s_1+s_2=2$.
Then
\begin{equation}\label{EQ:useful}
|x|^{-s_1}R_0(z)
|x|^{-s_2}\in \mathcal{B}(L^2(\Rn))
\end{equation}
for any $z\in \mathbb{C}\setminus [0,\infty)$.
\end{lem}
\begin{proof}
For completeness, we shall sketch the proof for $n=3$ and shall give the idea for $n >3.$
To begin with, we prepare Hardy inequality:
\begin{equation}\label{EQ:1-Hardy}
\left\| |x|^{-a} g\right\|_{L^2(\R^3)}\lesssim \left\| (-\Delta)^{\frac{a}{2}} g\right\|_{L^2(\R^3)}
\quad \text{for $0 \leq a<\frac32$,}
\end{equation}
and by duality,
\begin{equation}\label{EQ:2-Hardy}
\left\| (-\Delta)^{-\frac{a}{2}} h\right\|_{L^2(\R^3)}\lesssim \left\| |x|^a h\right\|_{L^2(\R^3)}
\quad \text{for $0 \leq a<\frac32$.}
\end{equation}
Noting that
\[
|R_0(z)g(x)|\le \frac{1}{4\pi}\int_{\R^3} \frac{|g(y)|}{|x-y|}\, dy=(-\Delta)^{-1}|g|(x)
\]
for any $z\in \mathbb{C}\setminus [0,\infty)$, we estimate
\begin{equation}\label{EQ:1-INE}
\begin{split}
\left\| |x|^{-s_1} R_0(z) |x|^{-s_2} f\right\|_{L^2(\R^3)}
= &\, \left\| |x|^{-s_1} \left|R_0(z) |x|^{-s_2} f \right| \, \right\|_{L^2(\R^3)}\\
\le &\, \left\| |x|^{-s_1} \left|(-\Delta)^{-1} |x|^{-s_2} |f| \right| \, \right\|_{L^2(\R^3)}
\end{split}
\end{equation}
for any $z\in \mathbb{C}\setminus [0,\infty)$.
Since
\[
\frac12<s_1<\frac32
\]
by our assumption on $s_1$ and $s_2$,
we can apply \eqref{EQ:1-Hardy} and \eqref{EQ:2-Hardy} to the right member of \eqref{EQ:1-INE}, and hence,
we estimate
\begin{equation}\label{EQ:2-INE}
\begin{split}
\left\| |x|^{-s_1} \left|(-\Delta)^{-1} |x|^{-s_2} |f| \right| \, \right\|_{L^2(\R^3)}
\lesssim&\,  \left\| (-\Delta)^{\frac{s_1}{2}} (-\Delta)^{-1} |x|^{-s_2} |f| \right\|_{L^2(\R^3)}\\
=&\, \left\| (-\Delta)^{-\frac{2-s_1}{2}} |x|^{-s_2} |f| \right\|_{L^2(\R^3)}\\
\lesssim &\, \left\|  |x|^{2-s_1-s_2} |f| \right\|_{L^2(\R^3)}\\
=&\, \|f\|_{L^2(\R^3)}.
\end{split}
\end{equation}
Combining \eqref{EQ:1-INE} and \eqref{EQ:2-INE}, we get the required estimate \eqref{EQ:useful}
for $n=3$.

For $n \geq 4$ one can represent $x=(x^\prime,x^{\prime\prime})$ where $x^\prime \in \mathbb{R}^{3}$
and one can use Plancherel's identity in $x^{\prime\prime}$ and then apply the estimate for $n=3.$
\end{proof}

In particular case for $s_1=s_2 = s >1$, we get
\begin{equation*} 
    \langle x \rangle^{-s}R_0(z)
\langle x \rangle^{-s}\in \mathcal{B}(L^2(\Rn))
\end{equation*}
for any $z\in \mathbb{C}\setminus [0,\infty)$.


\begin{thebibliography}{99}

\bibitem{Comech}
N. Boussaid and A. Comech, Nonlinear Dirac equation, Spectral stability of solitary waves,
Mathematical Surveys and Monographs, \textbf{244}. American Mathematical Society, Providence, RI, 2019.


%
\bibitem{Burq-Acta}
N.~Burq,
\textit{D\'ecroissance de l'\'energie locale de l'\'equation des ondes
pour le probl$\grave{e}$me ext\'erieur et absence de r\'esonance au voisinage
du r\'eel,} Acta Math. \textbf{180} (1998), 1--29.

%
\bibitem{Bu2004} N. Burq, \textit{Smoothing effect for Schr\"odinger boundary value problems,}
Duke Math. J. \textbf{123} (2004), no. 2, 403--427.


%
\bibitem{VoCa2002}
F. Cardoso and G. Vodev, \textit{Uniform estimates of the resolvent of the Laplace-Beltrami operator on infinite volume Riemannian
manifolds. {\rm II},}
Ann. Henri Poincar\'e \textbf{3}, (2002), no. 4, 673--691.

\bibitem{Duistermaat}
J. J. Duistermaat, \textit{Fourier integral operators.} Translated from Dutch notes of a course given at Nijmegen University, February 1970 to December 1971,
Courant Institute of Mathematical Sciences, New York University, New York, 1973.

\bibitem{Dunford}
 N. Dunford and J. T. Schwartz, Linear Operators I.
New York; Wiley-Interscience, 1966.

\bibitem{Evans} L. C. Evans, Partial Differential Equations. Graduate Studies in Mathematics, \textbf{19}.
American Mathematical Society, Providence, RI, 1998.

\bibitem{Georg-Matsuyama}
V. Georgiev and T. Matsuyama, \textit{Low frequency resolvent estimates for Dirichlet Laplacian on exterior domains},
AIP Conference Proceedings \textbf{2172}, 030012 (2019).

\bibitem{Georg-Visciglia}
V. Georgiev and N. Visciglia, \textit{Decay estimates for the wave equation with potential,}
Comm. Partial Differential Equations \textbf{28} (2003), 1325--1369.

\bibitem{Ginibre}
J. Ginibre and M. Moulin, \textit{Essential self-adjointness of many particle Schr\"odinger
Hamiltonians with singular two-body potentials,}
Ann. Inst. H. Poincar\'e Sect. A (N.S.) \textbf{21} (1974), 97--145.

\bibitem{Hormander2} L.~H\"ormander,
The Analysis of Linear Partial Differential Operators II. Springer-Verlag, Berlin, 1983.


\bibitem{Hormander} L.~H\"ormander,
The Analysis of Linear Partial Differential Operators III. Springer-Verlag, Berlin, 1985.

\bibitem{Ikawa} M. Ikawa,
\textit{Decay of solutions of the wave equation in the exterior of several convex bodies},
Ann. Inst. Fourier (Grenoble) \textbf{38} (1988), no. 2, 113--146.

\bibitem{Kerler}
C. Kerler,
\textit{Perturbations of the Laplacian with variable coefficients in exterior domains and differentiability of the resolvent},
Asymptot. Anal. \textbf{19} (1999), no. 3--4, 209--232.

\bibitem{LP}
P. D. Lax and R. S. Phillips, Scattering Theory. Academic Press, 2nd edition, 1990.


\bibitem{MM87}
R. Mazzeo and R. Melrose, \textit{Meromorphic extension of the resolvent on complete
spaces with asymptotically constant negative curvature,} J. Functional Analysis
\textbf{75} (1987), 260--310.


\bibitem{MS78}
R. B. Melrose and J. Sj\"ostrand,
\textit{Singularities of boundary value problems. {\rm I},}
Comm. Pure Appl. Math. \textbf{31} (1978), no. 5, 593--617.

\bibitem{MS82}
 R. B. Melrose and J. Sj\"ostrand,
\textit{Singularities of boundary value problems. {\rm II},} Comm. Pure Appl. Math.
\textbf{35}  (1982), no. 2, 129--168.



%
\bibitem{Mochizuki-2}K. Mochizuki,
Spectral and Scattering Theory for Second Order Elliptic
Differential Operators in an Exterior Domain.
Lecture Notes Univ. Utah, Winter and Spring 1972.

%
\bibitem{Mochizuki-1}
K. Mochizuki,
Spectral and scattering theory for second-order partial differential operators. Monographs and Research Notes in Mathematics.
CRC Press, Boca Raton, FL, 2017.


\bibitem{MRS}
C. Morawetz, J. Ralston and W. Strauss, \textit{Decay of solutions of the wave equation outside non-trapping obstacles,} Comm. Pure Appl.
Math. \textbf{30} (4) (1977), 447-- 508.

%

\bibitem{Moschidis}G. Moschidis,
\textit{Logarithmic local energy decay for scalar waves on a general class of asymptotically flat spacetimes,}
Ann. PDE \textbf{2} (2016), no.~1, Art. 5, 124 pp.


\bibitem{Petkov} V. Petkov, Scattering theory for hyperbolic operators, Elsevier Sci. Publ.,
North Holland, 1989.


\bibitem{R} J. Ralston,
\textit{Solutions of the wave equation
with localized energy,}
Comm. Pure App. Math. \textbf{22} (1969),
807--823.



\bibitem{TZ00} S.-H. Tang and M. Zworski, \textit{Resonance expansions of scattered waves,}
Comm. Pure Appl. Math. \textbf{53} (2000), 1305--1334.




\bibitem{Vodev2000}
G. Vodev, \textit{On the exponential bound of the cutoff resolvent,}
Serdica Math. J. \textbf{26} (2000), no. 1, 49--58.

\bibitem{Vodev2002}
G. Vodev, \textit{Local energy decay of solutions to the
wave equation for nontrapping metrics,}
Ark. Mat., \textbf{42} (2004), 379--397.





%
\bibitem{Wilcox}
C.~H.~Wilcox, \textit{Scattering theory for the d'Alembert equation
in exterior domains}, Lecture Notes Math. \textbf{442} (1975).





\end{thebibliography}
\end{document}